\documentclass[12pt]{article}
\usepackage[margin=1in]{geometry}
\usepackage{graphicx,url}
\usepackage{amssymb}
\usepackage{epstopdf}
\usepackage{amsthm}
\usepackage{amsmath}
\usepackage{multirow}
\usepackage{float}
\usepackage{color}
\usepackage{booktabs}
\usepackage[font=small,format=plain,labelfont=bf,up,textfont=it,up]{caption}
\usepackage{subfigure}
\usepackage[]{algorithm2e}

\newcommand{\R}{\mathbb{R}}

\DeclareMathOperator{\conv}{conv}
\DeclareMathOperator{\Proj}{Proj}

\title{Optimizing Damper Connectors for Adjacent Buildings}
\author{K. Bigdeli, W. Hare, J. Nutini and S. Tesfamariam}
\date{\today}

\begin{document}
\maketitle

\abstract

Many theoretical and experimental studies have used heuristic methods to investigate the dynamic behaviour of the passive coupling of adjacent structures. However, few papers have used optimization techniques with guaranteed convergence in order to increase the efficiency of the passive coupling of adjacent structures.
In this paper, the combined problem of optimal arrangement and mechanical properties of dampers placed between two adjacent buildings is considered. A new bi-level optimization approach is presented. The outer-loop of the approach optimizes damper configuration and is solved using the ``inserting dampers'' method, which was recently shown to be a very effective heuristic method. Under the assumption that the dampers have varying damper coefficients, the inner-loop finds the optimal damper coefficients by solving an $n$-dimensional optimization problem, where derivative information of the objective function is not available. Three different non-gradient methods are compared for solving the inner loop: a genetic algorithm (GA), the mesh adaptive direct search (MADS) algorithm, and the robust approximate gradient sampling (RAGS) algorithm. It is shown that by exploiting this new bi-level problem formulation, modern derivative free optimization techniques with guaranteed convergence (such as MADS and RAGS) can be used.  The results indicate a great increase in the efficiency of the retrofitting system, as well as the existence of a threshold on the number of dampers inserted with respect to the efficiency of the retrofitting system.
\\
\\

\noindent {\bf Keywords:} Seismic retrofitting, passive coupling, derivative-free optimization, bi-level optimization

\section{Introduction}

Increase in demand for residency and office buildings, coupled with limited land availability, has resulted in the construction of high-rise buildings in close proximity. During an earthquake, such closely spaced buildings are prone to pounding induced damages
\cite{tesfamariam2010, Ber86, KM97, cole2011}. By using damper connectors, the seismic vulnerability of adjacent buildings can be reduced, and thus, these pounding induced damages can be controlled \cite{azuma2006,bharti2010,xu1999}. Due to often having limited available resources, decision makers need to be able to optimize the number and placement of dampers.  In this paper, a new problem formulation and optimization approach are presented to solve the combined problem of finding the optimal arrangement and optimal mechanical properties of dampers placed between two adjacent buildings.

Considerable research has been reported on different damping devices, confirming the efficiency of damper connectors in mitigating the vibrations of structures \cite{xu2002, ying2003, bharti2010, qu2001, bhaskararao2006, ng2006, zhang1999, zhu2010, yang2003, zhang2000, ng2006, xu1999exp, yang2003}. The research on optimizing damper connections between adjacent buildings can be categorized into two categories: the placement of dampers; and the determination of mechanical properties for dampers \cite{spie_bigdeli2011, opt_loc_bigdeli2012}. The mechanical properties of dampers can be further subdivided into the kind of dampers used and the optimal {\em damper coefficients} required (see \cite{opt_loc_bigdeli2012, ok2008} or Section \ref{coefficients} herein for further details on damper coefficients). This paper focuses on the combined optimization of both the placement of dampers and the corresponding damper coefficients. This is set as a novel bi-level optimization problem of finding the optimal configuration of dampers and the corresponding optimal damper coefficients.  The configuration of dampers is a discrete optimization problem, whereas, the determination of optimal damper coefficients is a continuous optimization problem.

Several studies have presented results for optimizing damper coefficients \cite{xu1999, zhu2005, basili2007a, basili2007b, patel2010}. Xu et al.\ \cite{xu1999} considered multiple uniform dampers throughout the buildings, connecting every adjacent floor of the multiple degree-of-freedom (MDOF) structures. Zhu and Xu \cite{zhu2005} presented an analytical closed form solution for the damper coefficients of a fluid damper connecting two single degree-of-freedom (SDOF) structures. Basili and De Angelis \cite{basili2007a, basili2007b} studied optimal mechanical properties of nonlinear hysteretic dampers connecting SDOF and MODF structures. They presented explicit equations relating dissipation energy, relative displacement and relative acceleration of the mechanical properties of the dampers. Patel and Jangid  \cite{patel2010} assumed that the optimal damper coefficients were functions of the relative velocity between the structures; the damper coefficients increased from a small value for the base floor to a large value for the top floor.

The discrete optimization problem of damper positioning has also been the topic of several recent studies \cite{yang2003, ok2008, opt_loc_bigdeli2012}. For example, Yang and Lu \cite{yang2003} constructed a series of experiments to show that one can eliminate half of the dampers without compromising efficiency. However, a method to determine the optimal arrangement of the dampers was not included. Bigdeli et al.\ \cite{opt_loc_bigdeli2012} presented an optimization algorithm to find the optimal placement of dampers with a limited number of available dampers. All dampers were assumed to be identical with the same mechanical properties. The results of Bigdeli et al.\ \cite{opt_loc_bigdeli2012} also demonstrated that if all damper coefficients are assumed to be equal, then increasing the number of dampers does not necessarily increase the dynamic stability of a structure.  Moreover, under these conditions, increasing the number of dampers may actually exacerbate the dynamic behaviour of the buildings. This is in agreement with other studies \cite{ok2008}.

To the authors' knowledge, the only research that studies the combined damper location and coefficient selection problem is Ok et al. \cite{ok2008}. Ok et al.\ \cite{ok2008} examined a multi-objective optimization method using genetic algorithm (GA) for a set of coupled MDOF structures connected to each other by magneto-rheological (MR) dampers. The number of dampers and the voltage for the MR dampers installed at each floor were assumed as design parameters. Since a bounded domain for voltage was assumed, they allowed each floor to have more than one damper if the result was more effective.
Unlike this paper, Ok et al.\ used a single-level optimization framework, and thus, used a heuristic (specifically genetic algorithm) to seek optimality.  In this paper, a novel bi-level framework is considered to optimize the arrangement and mechanical properties of dampers placed between two adjacent buildings.  The objective function is set to {\em minimize} the {\em maximum inter-story drift} over all possible damper configurations. This is a common objective in seismic retrofitting literature \cite{LavanCimellaroReinhorn2008, opt_loc_bigdeli2012, Kanno2013}, however, it should be noted that other objective functions have been used \cite{yang2003, ok2008, GrecoLucchiniMarano2014}.  The optimization techniques and bi-level framework used herein are independent of the objective function, so could be readily applied to optimize other aspects of seismic retrofitting.  Such exploration is left for future research.

The bi-level approach presented in this paper uses an outer-loop that seeks to optimize the combinatorial problem of where to place dampers, and an inner-loop that seeks to solve the continuous problem of determining the optimal design for each damper.  This bi-level approach presents one key advantage over a single level approach, namely, the inner-loop can be solved using modern derivative-free optimization (DFO) software, and therefore some assurance of (local) optimality can be attained.  The approach is tested on a suite of 150 test-cases.  Results support the effectiveness of the approach and demonstrate the importance of using high quality DFO software (opposed to heuristic methods).  Results further demonstrate that there exists a threshold on the number of dampers inserted with respect to the efficiency of the retrofitting system.  That is, maximal efficiency can be achieved using a limited number of dampers, provided that the internal damper designs are fully optimized.

The remainder of this paper is organized as follows. In Section \ref{model}, the physical model and the optimization problem are presented. In Section \ref{Algorithms}, a discussion of each level of the bi-level problem is presented, as are the details of the optimization methods used in this paper. In Section \ref{Numerics}, numerical results for 150 test problems are presented. Some conclusions are provided in Section \ref{Conclusions}.

\section{Physical Model \& Optimization}\label{model}

\begin{figure}[ht]
\centering
\includegraphics[height=15.1cm,width=7cm]{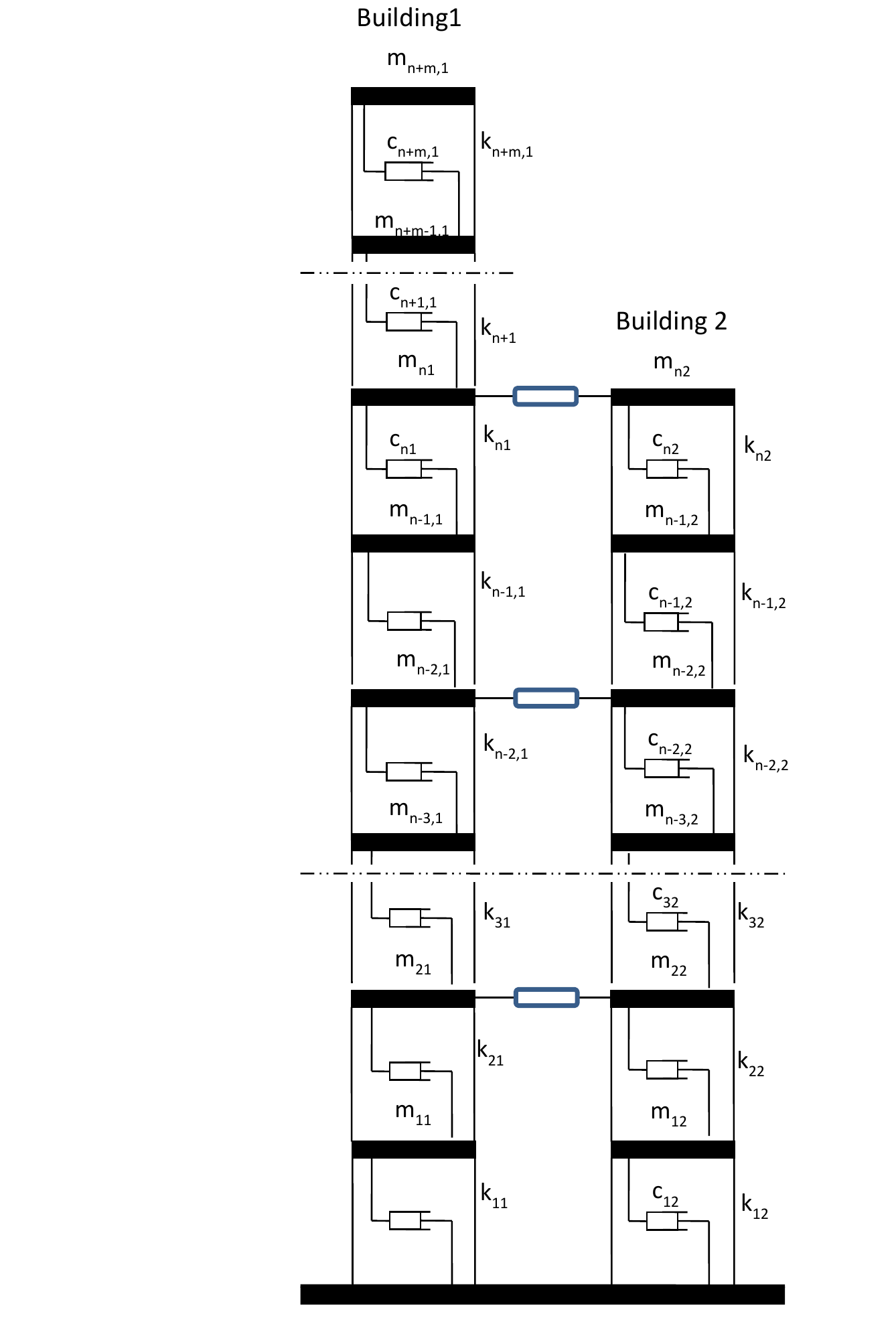}
\caption[Building model]{Model of adjacent buildings with damper connectors.}
\label{figure:1}
\end{figure}

Figure \ref{figure:1} provides an illustrative model of two adjacent buildings connected by dampers. In this paper, buildings are modelled as 1-dimensional systems assuming the center of rigidity and the center of mass of adjacent floors are in the same plane, and the viscous dampers (with variable damper coefficients) are connected at the floor level.  Each building is modelled as a multi degree-of-freedom (MDOF) system consisting of lumped masses (the mass of each floor), linear springs (stiffness of the columns), and linear viscous dampers \cite{xu1999,opt_loc_bigdeli2012}.  The ground motion and dynamic response of the buildings are assumed to be unidirectional. Both buildings are assumed to be symmetric in plane (i.e., their centers of mass located in the same plane), and as a result, the effect of torsional vibrations is not considered.  Equal building heights are not required, and consequently, the maximum number of dampers that can be placed is controlled by the shortest building.

Buildings $1$ and $2$ have $n + m$ and $n$ stories and are connected by $n_d$ dampers.  Thus, the dynamic model for both structures is  a $2n + m$ degree-of-freedom system. Let $X(t) \in \R^N$, where $N=2n+m$, be the vector of displacements of each floor at time $t$.  The governing equation of the system can be expressed as
\begin{equation}\label{eq:1}
M \ddot X(t)+(C+C_d) \dot X(t)+ K  X(t)= M  E \,  g(t),
\end{equation}
where matrices $M \in \R^{N\times N}$, $C \in \R^{N\times N}$ and $K \in \R^{N\times N}$ are generated by the given mass, damping, and stiffness factors of the buildings, respectively.  The vector $E \in \R^N$ is a vector of ones, and the function $g: \R \rightarrow \R$ is the ground acceleration during the earthquake.  The matrix $C_d  \in \R^{N\times N}$ is constructed using the damper coefficients and the locations of the $n_d$ dampers.  Let $c_d \in \R^n$, $c_d \geq 0$, be the vector of damper coefficients for each floor, where any floor $i$ that is without a damper has coefficient $c_{d,i}=0$.  The matrix $C_d$ takes the form
\begin{equation}
	C_d = \left[ \begin{array}{cccccccc}
	\mathrm{diag}(c_d) & \mathrm{zero}(n,m) & -\mathrm{diag}(c_d) \\\\
	\mathrm{zero}(m,n) & \mathrm{zero}(m,m) & \mathrm{zero}(m,n) \\ \\
	-\mathrm{diag}(c_d) & \mathrm{zero}(n,m) & \mathrm{diag}(c_d)
	\end{array}\right],
\end{equation}
where $\mathrm{diag}(x)$ is the diagonal matrix whose entries coincide with the vector $x$, and $\mathrm{zero}(a,b)$ is an $a \times b$ zero matrix.

Note that equation \eqref{eq:1} describes the motion in the time domain. Considering the spectral density of the ground excitation, the equation of motion can be written in the frequency domain as
\begin{equation}\label{eq:2}
e^{i\omega t} \left[ -M \omega^2 X(\omega)+(C+C_d)i \omega X(\omega) +K X(\omega) \right]=-M E \sqrt{S_g(\omega)}e^{i\omega t},
\end{equation}
where the response of the building is given by
\begin{equation}\label{eq:3}
X(\omega)=\left[ -M \omega^2 +(C+C_d)i \omega + K  \right]^{-1} \times \left[-M E \sqrt{S_g(\omega)} \right].
\end{equation}

In this paper, a Kanai-Tajimi filtered white noise function is used for the spectral density function of ground acceleration:
\begin{equation}\label{eq:4}
S_g(\omega)=\frac{ 1+4\zeta_g^2  \left(\frac{\omega}{\omega_g}\right)^2  }  { \left(1-\left( \frac {\omega}{\omega_g}\right)^2 \right)^2   +4\zeta_g^2  \left(\frac
{\omega}{\omega_g}\right)^2} S_0,
\end{equation}
where $\omega_g$, $\zeta_g$ and $S_0$ represent dynamics characteristics and the intensity of the earthquake. These parameters are chosen based on geological characteristics.

For a given vector of damper coefficients, $c_d$, a numerical approximation of the standard deviation of the displacement response for the $i^{th}$ floor of building $b$ is possible, and is given by
\begin{equation}\label{eq:5}
\sigma_{ib} = \left[ \int_{-\infty}^{+\infty} \|x_{ib}(\omega)\|^2 d\omega \right]^{\frac{1}{2}},
\end{equation}
where $x_{ib}(\omega)$ is the component of $X(\omega)$ corresponding to floor $i$ of building $b$.  The value $\sigma_{ib}$ is used to calculate the {\em inter-story drift} for each floor, which is denoted by $f_{ib} = \big (\sigma_{ib} - \sigma _{(i-1)b} \big )^2$, where $\sigma_{0b}$ is defined as $0$.  This in turn defines the maximum inter-story drift as
\begin{equation}\label{eq:6}
F = \max \bigg \{ \max_{i=1,\dots,n+m} \big \{ (\sigma_{i1} - \sigma _{(i-1)1})^2 \big \} \: , \:  \max_{i=1,\dots,n} \big\{ (\sigma_{i2} - \sigma _{(i-1)2} )^2 \big\} \bigg \}.
\end{equation}
The objective is to determine the optimal configuration of dampers and damper coefficients in order to {\em minimize} the {\em maximum inter-story drift}.

As an analytic solution to equation \eqref{eq:5} is unavailable, a numerical approximation of $\sigma_{ib}$ is required. To do this, upper and lower limits of $\pm 20 {rad/s}$ are imposed, and a trapezoidal rule approximation is applied to the integral in equation \eqref{eq:5} with a step size of $0.02$.  (Previous studies show that the effect of frequencies greater than $20 {rad/s}$ on the response of the structure is negligible \cite{xu1999}).

In summary, to place $n_d$ damper connectors between two buildings of heights $n+m$ and $n$, an optimization problem of the following form is considered:
	\begin{equation} \min_{c_d \in \R^n_+} F({c_{d}})  ~~ s.t. ~~
			  c_{d,j} = 0 ~~\mbox{for at least $n-n_d$ values of}~j,
	 \end{equation}
where
	\begin{equation}
	F(c_{d}) = \max \big \{ f_{ib}(c_d) : i = 1,\dots, n +m, b=1; i= 1, \dots, n, b = 2 \big \}
	\end{equation}
and each $f_{ib}$ is numerically approximated using computer simulation. The bi-level optimization approach presented in the next section is designed for this problem formulation: the outer-loop seeks to optimize the combinatorial problem of which values of $c_{d,j}$ are $0$, and the inner-loop seeks to solve the continuous problem of determining the optimal values of the damper coefficients.

\section{A Bi-level Approach}\label{Algorithms}

\begin{figure}[ht]
\centering
\includegraphics[scale=0.75]{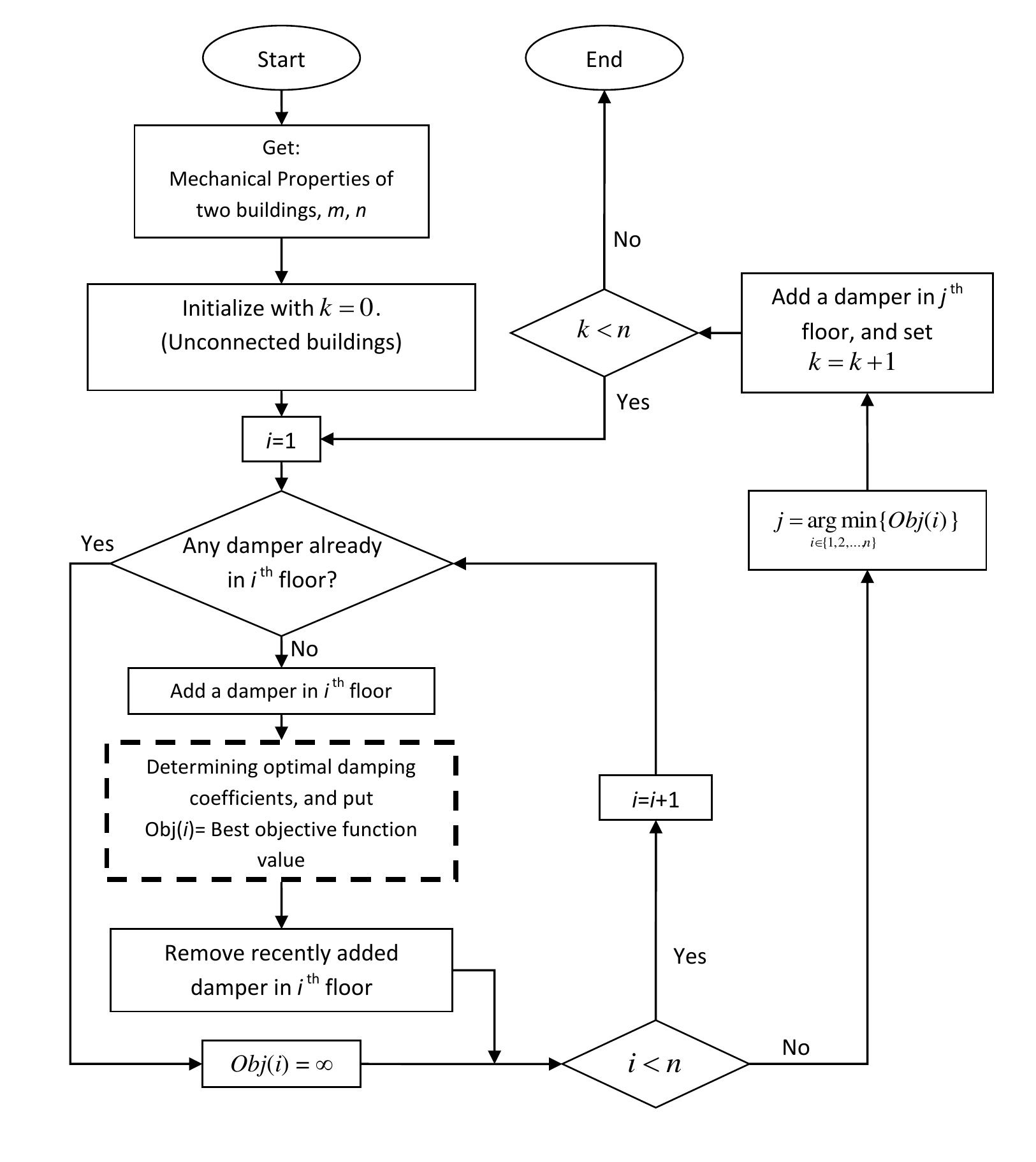}
\caption{Schematic of the bi-level optimization problem: outer-loop in full, inner-loop (dashed box) in brief.}
\label{figure:schematic}
\end{figure}

As discussed in the previous section, the optimization problem in this paper is considered as a bi-level optimization problem with an inner continuous optimization algorithm and an outer discrete optimization algorithm. The inner-loop uses a non-gradient based method to find an optimal set of damper coefficients for a fixed configuration and a fixed number of dampers. Based on the research in \cite{opt_loc_bigdeli2012}, the outer-loop uses a heuristic optimization algorithm, which seeks the optimal configuration of dampers.   A schematic outline of the presented bi-level approach is given in Figure \ref{figure:schematic}.

Beginning at initialization, the algorithm obtains the mechanical properties of the two adjacent buildings under consideration. The buildings are initially considered as two unconnected structures, and the iteration count $k$ is set to $0$. From here, the algorithm initializes the first iteration of the outer-loop, which optimizes damper location.

The heuristic algorithm used for the outer-loop (as schematically shown in full in Figure \ref{figure:schematic}) is the inserting dampers method from \cite{opt_loc_bigdeli2012}. In \cite{opt_loc_bigdeli2012}, this method was shown to be the most effective among others at finding the optimal configuration of dampers.  In addition, this approach has the advantage that the optimization of $n_d$ dampers automatically provides the solutions to the optimization of $1, 2, ...,$ and $n_d - 1$ dampers.

During the first iteration, the algorithm checks all possible locations, i.e., all floors $i= 1, \dots n$, to put the first damper by determining the optimal damper coefficients (inner-loop, dashed line in Figure \ref{figure:schematic}) for each location. Following each inner-loop iteration, the resulting objective value (minimum over the maximum inter-story drifts for the current fixed damper configuration) is compared to the best objective value for the current floor. In the first iteration, this step consists of initializing the objective value for each floor. For subsequent iterations, if placing the damper on the current floor in the current configuration produces a smaller objective function value than the current best, then the stored best objective function value is updated accordingly for that floor. The damper is then removed from the current floor, and a damper is placed on the next (consecutive) floor. This continues until the top floor is reached.

Once all the floors have been cycled through (condition $i < n$ is not satisfied), the algorithm computes an element of the argmin over the resulting set of $n$ objective function values.  (Recall, the argmin, or the argument of the minimum, is the set of all minimizers for an optimization problem.) The solution to this problem, $j$, indicates that the placement of a damper on floor $j$ results in the minimal objective value for the current overall damper configuration. Thus, a `permanent' damper is inserted at floor $j$. From here, the iteration count is increased and this iterative procedure is repeated, starting again at floor 1, until all available dampers are inserted into the structure. (In order to aid the inner-loop, when possible the solution to the previous iteration is used as a warm-start for the solver used in the inner loop (see Section \ref{coefficients})).  If a `permanent' damper has been inserted in a previous iteration at the current floor, then when that floor is reached in the cycle, the algorithm assigns an objective function value of infinity for this floor. It then proceeds onto the next step: either continuing onto the next floor, or computing the argmin over the objective values. This ensures that any floor with a `permanent' damper will not be selected again for damper insertion when calculating the argmin over the objective function values. This iterative process continues until the maximum number of allowable dampers have been placed.


In the next section, details of the methods used to solve for the optimal damper coefficients in the inner-loop are presented.

\subsection{Damper Coefficients Optimization}\label{coefficients}

The main purpose of the inner-loop of the optimization algorithm (the dashed box in Figure \ref{figure:schematic}) is to find optimal damper coefficients for a fixed configuration of damper connectors. As derivative information of the objective function is not available for the optimization problem, the inner-loop requires the use of a non-gradient based method.

Many non-gradient optimization options exist.  These can be widely split into two categories, heuristic optimization methods and derivative-free optimization (DFO) methods.  Here, DFO refers to methods that are mathematically derived and studied to provide (theoretical) proof of convergence to (local) minimizers, whereas heuristic methods are any other non-gradient based methods that do not fit this definition.   For a thorough introduction into several well-known DFO frameworks, see \cite{hare2013, CSV09}.

Due to their versatility, heuristic optimization methods are widely used in structural engineering. One such method is the genetic algorithm (GA) \cite{ok2008}. However, heuristic methods do not guarantee convergence to (locally) optimal solutions. As such, there has been a recent increase in the use of derivative-free optimization techniques that guarantee optimality. In this work, the mesh adaptive direct search (MADS) algorithm \cite{audet2006} is examined, as implemented in MATLAB's global optimization toolbox; as well as a novel robust approximate gradient sampling (RAGS) algorithm \cite{hare2013derivative} that is specifically designed for finite minimax problems \cite{hare2013derivative}. The stochastic based GA is used as a baseline comparison to previous studies that have principally employed this method.  A detailed description of each algorithm is given in Section \ref{Algorithms}.  Greater detail on each of these methods is provided next.

\subsubsection{Genetic Algorithm}
The genetic algorithm (GA) is a popular heuristic search method. It has been argued to be a reasonably efficient method, particularly in engineering applications (\cite{alkhatib2004, singh2002, hadi1998, ok2008} and references therein). A simple GA consists of several steps, including the generation of initial points, selection, competition and reproduction \cite{vose1999}. A brief description of the genetic algorithm follows.
\bigskip \\
\begin{tabular}{ll}
{\bf procedure} GeneticAlgorithm\\
{\bf begin} \\
\hspace{5mm}Initialize and evaluate random population $P(t)$; \\
\hspace{5mm}{\bf while} {\em stopping conditions not satisfied} {\bf do} \\
\hspace{5mm}{\bf begin} \\
\hspace{5mm}\hspace{5mm}Mutate and crossover $P(t)$ to yield $C(t)$; \\
\hspace{5mm}\hspace{5mm}Evaluate $C(t)$; \\
\hspace{5mm}\hspace{5mm}Select $P(t+1)$ from $C(t)$ and elite individuals from $P(t)$; \\
\hspace{5mm} {\bf end} \\
{\bf end}
\end{tabular}
\bigskip \\
In this paper, the GA is used in a standard form that is included in the MATLAB global optimization toolbox \cite{MATLAB}.  All parameter values and settings are the MATLAB default choices, future research may explore alternate parameter selections.

\subsubsection{Mesh adaptive direct search}
The mesh adaptive direct search (MADS) method \cite{audet2006} is a sub-category of pattern search (PS) methods. A brief description of a general pattern search method follows.
\bigskip \\
\begin{tabular}{ll}
{\bf procedure} PatternSearch\\
{\bf begin}  \\
\hspace{5mm}Initialize $x_0$ and a set of directions $\mathcal{D}$; \\
\hspace{5mm}{\bf while} {\em stopping conditions not satisfied} {\bf do} \\
\hspace{5mm}{\bf begin} \\
\hspace{5mm}\hspace{5mm}{\bf Search} for a point with $f(x) < f(x_k)$ (optional);\\
\hspace{5mm}\hspace{5mm}{\bf Poll} points from $\{x_k + \alpha_k d : d \in D_k (\in \mathcal{D}) \}$;\\
\hspace{5mm}\hspace{5mm}{\bf if} $f(x_k + \alpha_k d_k ) < f(x_k)$\\
\hspace{5mm}\hspace{5mm}\hspace{5mm}Stop polling;\\
\hspace{5mm}\hspace{5mm}\hspace{5mm}$x_{k+1} \leftarrow x_k + \alpha_k d_k$;\\
\hspace{5mm}\hspace{5mm}{\bf else}\\
\hspace{5mm}\hspace{5mm}\hspace{5mm} $x_{k+1} \leftarrow x_k$;\\
\hspace{5mm}\hspace{5mm}\hspace{5mm}{\bf Update} mesh parameter $\alpha_k$;\\
\hspace{5mm}\hspace{5mm}{\bf end}\\
\hspace{5mm}{\bf end} \\
{\bf end}
\end{tabular}
\bigskip\\
Specific to MADS, randomly rotated bases are used in each iteration to provide a more robust convergence. In this paper, a MATLAB interface is employed with the version of MADS that is implemented in the NOMAD project \cite{NOMAD}.  All parameter values and settings are the MADS default choices, future research may explore alternate parameter selections.

\subsubsection{Robust approximate gradient sampling}
The robust approximate gradient sampling algorithm (RAGS algorithm) is a derivative-free optimization algorithm that exploits the smooth substructure of the finite minimax problem,
\[
	\min_x F(x) \text{ where } F(x)=\max\{ f_i : i = 1, \dots , N \}.
\]
The general concept of the RAGS algorithm relies on the definition of the {\em active set} of a finite max function $f$ at a point $\bar{x}$,
\[
	A(\bar{x}) = \{ i : F(\bar{x}) = f_i (\bar{x}) \}.
\]
Loosely speaking, the {\em subdifferential} of $F$ at a point $\bar{x}$ is the set of all possible gradients. For a finite max function, the subdifferential, as shown in \cite{Clarke90}, is given by
\begin{equation}\label{sub}
	\partial f(\bar{x}) = \conv\{\nabla f_i (\bar{x})\}_{i \in A(\bar{x})},
\end{equation}
and the direction of steepest descent can be defined via $-\Proj(0|\partial f(\bar{x}))$. Although the direction of steepest descent is fine, it tends to get stuck on non-differentiable ridges of the function.

In 2005, Burke et al.\ \cite{BLO05} introduced a robust gradient sampling algorithm. This algorithm uses information from around the current iterate to help minimize along non-differentiable ridges of nonsmooth functions. A brief description of the RAGS algorithm follows.
\bigskip \\
\begin{tabular}{ll}
{\bf procedure} RAGS\\
{\bf begin}\\
\hspace{5mm}Initialize $x_0$, search radius $\Delta_0$, Armijo-like parameter $\eta$ and other parameters; \\
\hspace{5mm}{\bf begin} \\
\hspace{10mm}Generate a set of $n+1$ points;\\
\hspace{10mm}Use points to generate robust approximate subdifferential $G_Y^k$;\\
\hspace{10mm}Set search direction $d_Y^k = \Proj(0|G_Y^k)$;\\
\hspace{10mm}{\bf if} $\Delta_k$ small, but $|d^k|$ large\\
\hspace{15mm}Carry out line search: find $t_k>0$ such that $f(x^k +t_k d^k) < f(x^k) - \eta t_k |d^k|^2$;\\
\hspace{20mm}{\bf Success}: update $x^k$ and loop;\\
\hspace{20mm}{\bf Failure}: decrease accuracy measure and loop;\\
\hspace{10mm}{\bf else if} $\Delta_k$ large\\
\hspace{15mm}Decrease $\Delta_k$ and loop;\\
\hspace{10mm}{\bf else}\\
\hspace{15mm}Terminate;\\
\hspace{10mm}{\bf end}\\
\hspace{5mm}{\bf end}\\
{\bf end} \\
\end{tabular}
\bigskip\\
The RAGS algorithm uses approximate gradients to adapt the robust gradient sampling algorithm to a DFO setting.  Like MADS, the RAGS algorithm is proven to converge to a local minimizer. Readers are referred to \cite{hare2013derivative} for further information.  All parameter values and settings are the RAGS default choices, future research may explore alternate parameter selections.

\subsection{Warm-starting}

All of the algorithms tested allow the user to input an initial point. Thus, once the outer-loop has completed at least one iteration, the solution to the past outer-loop configuration can be used to `warm-start' the inner-loop computation.  Specifically, the damper coefficients are initialized by setting $c_{d,i}^j = c_{d,i}^{j-1}$ if there was a damper in position $i$ during outer-loop $j-1$, and $c_{d,i}^j = 0$ if there was not a damper in position $i$ during outer-loop $j-1$.  Two versions of each algorithm are considered, for six algorithms total: GA using a random initial point for each new inner-loop (denoted GA$_r$), GA using the warm-start initial point for each new inner-loop (denoted GA$_w$), MADS using a random initial point (denoted MADS$_r$), MADS using the warm-start initial point (denoted MADS$_w$), RAGS using a random initial point (denoted RAGS$_r$), and RAGS using the warm-start initial point (denoted RAGS$_w$).

\section{Numerical Results}\label{Numerics}

In this section, a summary of results are presented for various numerical problems for damper-connected structures. The solution times and quality measures of the previously presented non-gradient based methods are compared.

\subsection{Test Problems}

In order to compare the presented methods, three different sets of mechanical properties and five different sets of heights for the two buildings are considered.  It should be noted that the mechanical properties are taken from previous studies in the field and are reasonable examples of buildings requiring seismic retrofitting \cite{xu1999} \cite{opt_loc_bigdeli2012}.  However, the building heights are artificial, and selected to be representative of a variety of scenarios.  These are provided in Tables \ref{table:1} and \ref{table:2}.

\begin{table}[ht]
  \centering
  \caption{Mechanical Properties}
    \begin{tabular}{lccccccccc}
    \toprule
          & \multicolumn{3}{c}{Building $a$} &       & \multicolumn{3}{c}{Building $b$} & \\
    \cmidrule(r){2-9}
          & ma (kg) & k­a (N/m) & ca (N.s/m) & \textit{} & mb (kg) & k­b (N/m) & cb (N.s/m) \\
    \cmidrule(r){2-9}
    Set I & 1.29E+06 & 4.00E+09 & 1.00E+05 &       & 1.29E+06 & 2.00E+09 & 1.00E+05 \\
    Set II & 2.60E+06 & 1.20E+10 & 2.40E+06 &       & 1.60E+06 & 1.20E+10 & 2.40E+06 \\
    Set III & 4.80E+06 & 1.60E+10 & 1.20E+06 &       & 4.00E+06 & 2.30E+10 & 1.20E+06 \\
    \bottomrule
    \end{tabular}%
  \label{table:1}%
\end{table}%

\begin{table}[ht]
  \centering
  \caption{Building Heights}
    \begin{tabular}{crcr}
    \toprule
   Case  & $f_a$    & $f_b$    &  \\
    \midrule
    1 & 10    & 10    &   \\
    2 & 10    & 20    &   \\
    3 & 20    & 10    &   \\
    4 & 10    & 40    &   \\
    5 & 40    & 10    &   \\
    \bottomrule
    \end{tabular}%
  \label{table:2}%
\end{table}%

In Table \ref{table:2}, $f_a$ and $f_b$ represent the number of floors for buildings $a$ and $b$, respectively. For all numerical examples, to generate the ground excitation spectrum, the following values are used for the ground acceleration parameters in equation \eqref{eq:4}: $\omega_g = 15 ~rad/s$, $\zeta_g=0.6$, $\omega_k=1.5 ~rad/s$ and $S_0 = 4.65\times10^{-4} ~m^2 / rad.s^3$. (These parameter values are the same as those used in \cite{opt_loc_bigdeli2012} and \cite{xu1999}.) For each of the 3 sets of mechanical properties, the number of dampers changes from $1$ to $10$. Therefore, incorporating all 5 building height combinations, a total of 150 test problems are generated, representing a wide range of situations. Each problem is solved via a combination of the inserting dampers method and a non-gradient based method (either GA, MADS or RAGS). Optimal arrangements and damper coefficients, as well as corresponding objective function values are determined.

As an example output, for Building Height 1, Material Set I, using 4 dampers, the optimal configurations represented in vector form, as solved via GA, MADS, and RAGS were as follows:
    $$\begin{array}{rll}
    \mbox{GA:} & [\mathbf{0}, \mathbf{0}, 2.4331, 0.4821, 1.5187, \mathbf{0}, \mathbf{0}, \mathbf{0}, \mathbf{0}, 0.2146] * 10^{7}, \\
    \mbox{MADS:} & [\mathbf{0}, \mathbf{0}, 2.4179, 0.1000, 1.6550, \mathbf{0}, \mathbf{0}, \mathbf{0}, \mathbf{0}, 0.2257] * 10^{7}, \\
    \mbox{RAGS:} & [\mathbf{0}, \mathbf{0}, 2.4188, 0.1135, 1.6505, \mathbf{0}, \mathbf{0}, \mathbf{0}, \mathbf{0}, 0.2099] * 10^{7}, 
    \end{array}$$
where a bold zero ($\mathbf{0}$) denotes a floor that has no damper.  Notice that, while the damper coefficients differ, they are similar, and all methods resulted in dampers on floors 3, 4, 5, and 10.  Also note that the damper coefficients differ from floor to floor, emphasizing the need for multi-variable optimization.

\subsection{Solution Time and Quality}\label{ss:time}

Tables \ref{table:dataMatIrandom} to \ref{table:dataMatIIIwarm} in Appendix A show the number of function calls required and the optimal objective function values obtained using various methods. For the sake of brevity, optimal design variables, including configurations of dampers and damper coefficients, are not included. The full data is available upon request by contacting the corresponding author.

Note that in Tables \ref{table:dataMatIrandom} to \ref{table:dataMatIIIwarm}, instead of reporting actual solution times in seconds, the number of performed simulations is reported. It is worth noting that each simulation takes approximately 2 seconds, regardless of the details or dimension of the problem.

\begin{figure}[ht]
\centering
\subfigure[Objective value at each function call.]{
\includegraphics[width=2.5in]{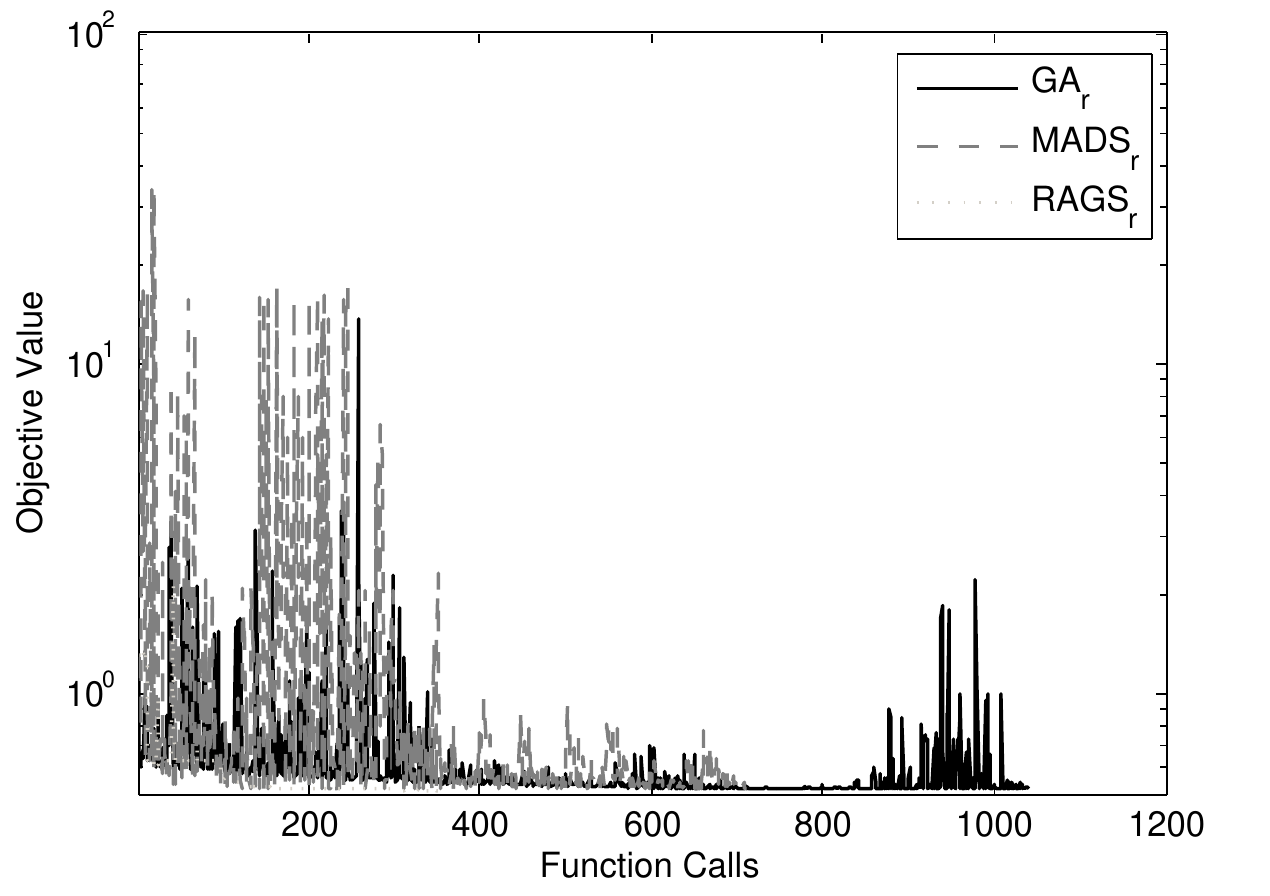}
\label{figure:randconvergenceratea}
}
\subfigure[Minimum objective value at each function call.]{
\includegraphics[width=2.5in]{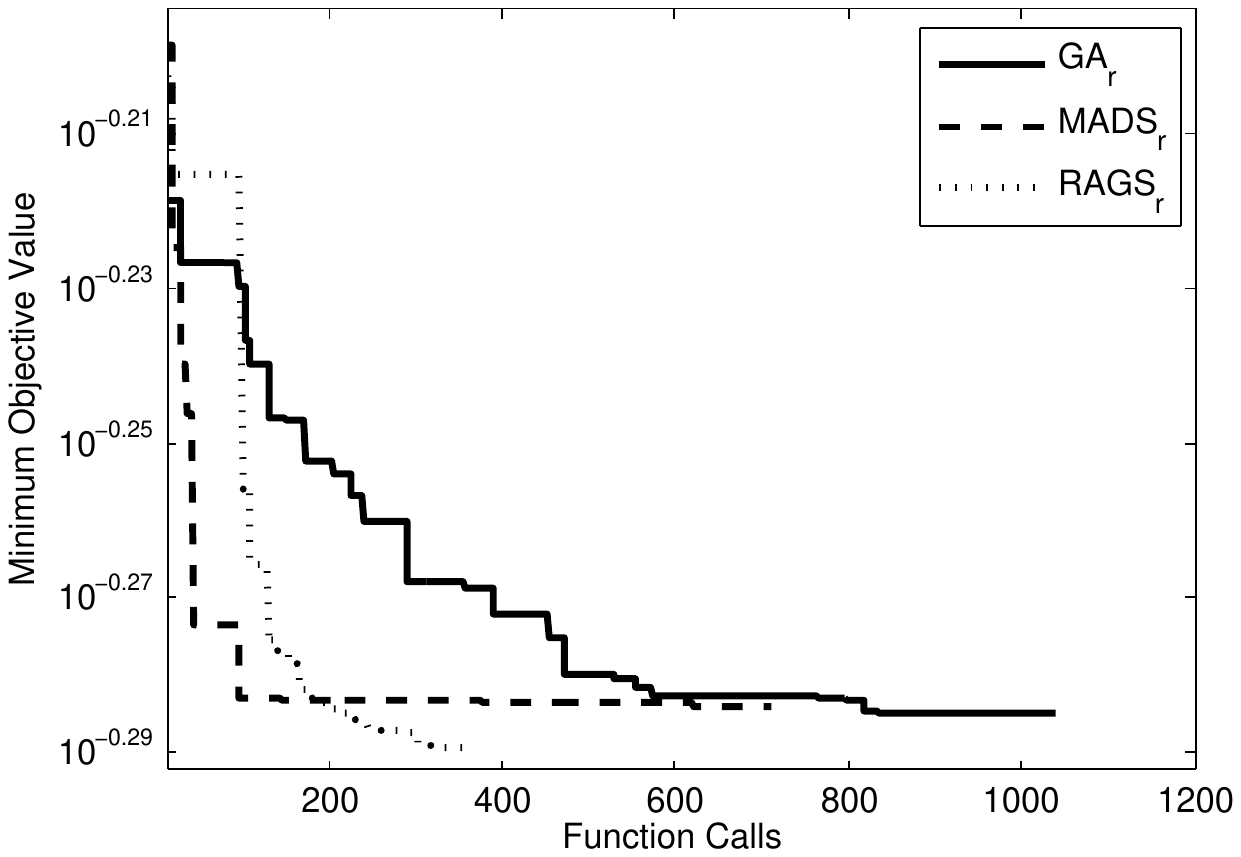}
\label{figure:randconvergencerateb}
}
\caption{Objective values for Material Set I and Building Heights 1 when algorithms use random initial points.}
\label{Figure:randconvergencerate}
\end{figure}

\begin{figure}[ht]
\centering
\subfigure[Objective value at each function call.]{
\includegraphics[width=2.5in]{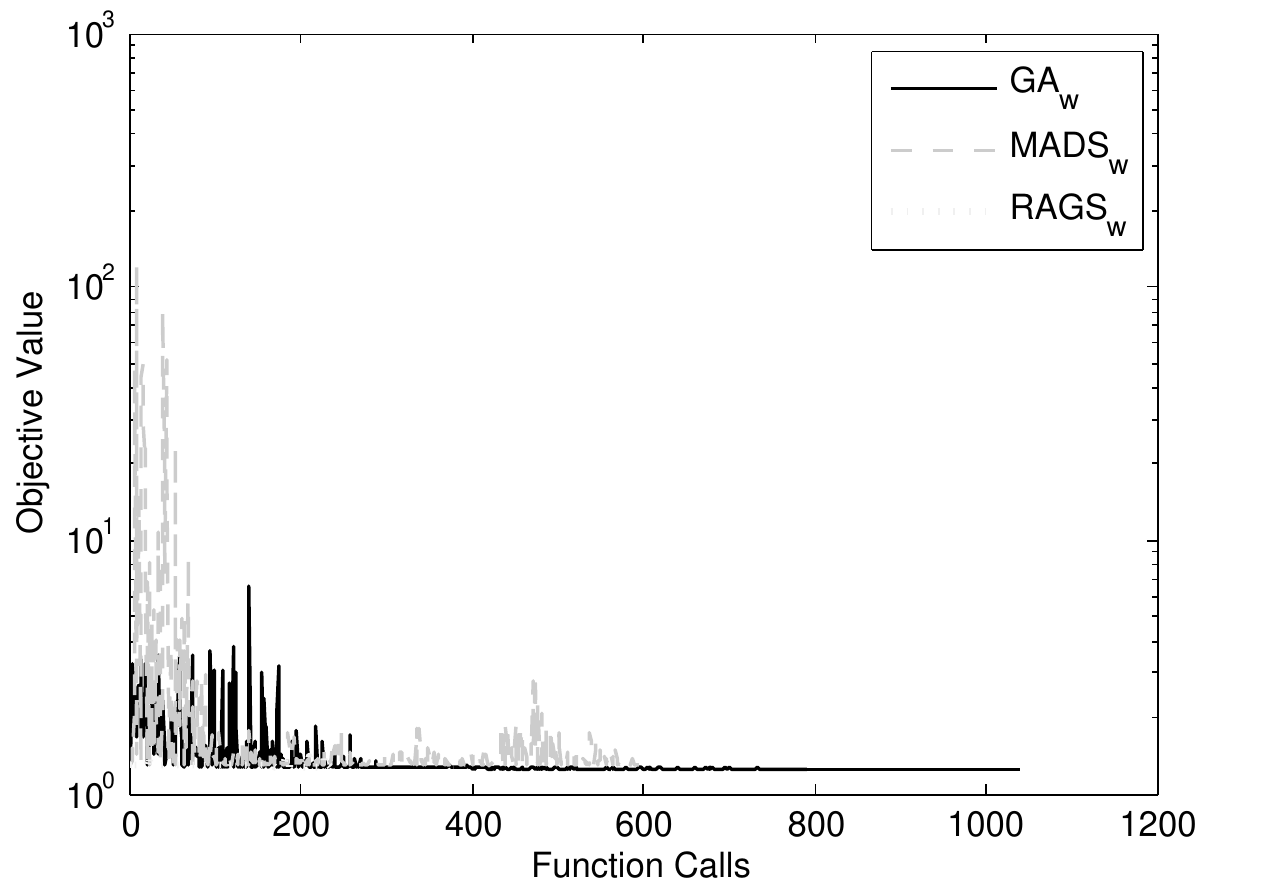}
\label{figure:warmconvergenceratea}
}
\subfigure[Minimum objective value at each function call.]{
\includegraphics[width=2.5in]{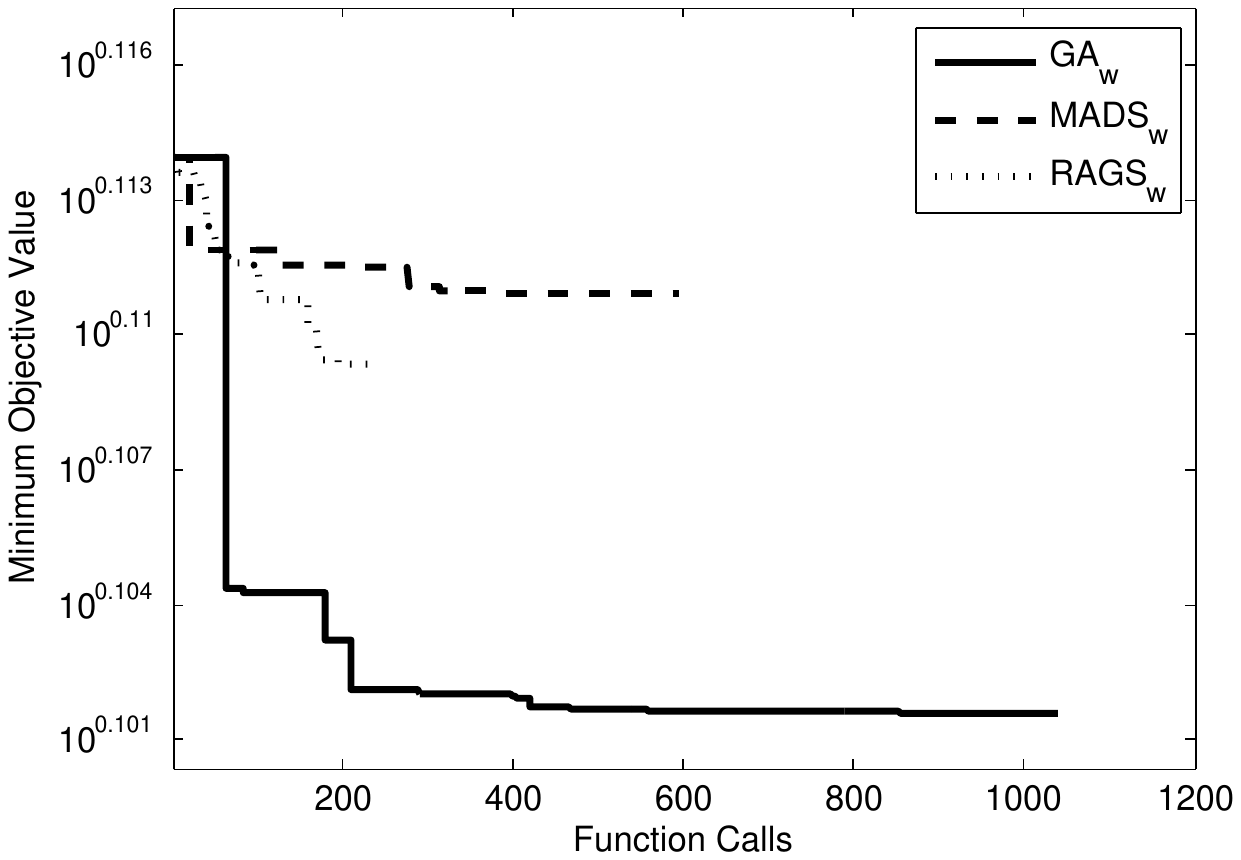}
\label{figure:warmconvergencerateb}
}
\caption{Objective values for Material Set I and Building Heights 1 when algorithms use warm-start initial points.}
\label{Figure:warmconvergencerate}
\end{figure}

Clearly, the rate of convergence is a crucial factor for any optimization algorithm. As a first comparison of convergence rates, in Figures \ref{Figure:randconvergencerate} and \ref{Figure:warmconvergencerate}, for Building Height 1, Material Set I, and the case when all adjacent floors are connected, the objective value (a) and the minimum objective value (b) for each function call are plotted. Figure \ref{Figure:randconvergencerate} displays the three algorithms using random initial points and Figure \ref{Figure:warmconvergencerate} displays the three algorithms using warm-start initial points.  For brevity, other figures displaying similar results are omitted.

In Figure \ref{figure:randconvergenceratea}, it can be seen that, GA is a stochastic based method.  In particular, it evaluates the objective function at a wide range of points, resulting in a large range of function values, even after convergence is essentially established.  For the MADS algorithm, a similar variation in objective value range is seen, with multiple spikes in the objective value as the number of function calls increases. Looking closely, it can be seen that RAGS converges with minimal objective value variation. This is because RAGS lacks any global search heuristics.

Examining Figures \ref{Figure:randconvergencerate} and \ref{Figure:warmconvergencerate}, one notes that all three algorithms do fairly well on these problems.  In Figure \ref{Figure:randconvergencerate}, RAGS outperforms the other methods, while in Figure \ref{Figure:warmconvergencerate}, GA does extremely well.  However, it should be noted that these figures only represent 1 out of 150 test problems.  In order to investigate the overall performance of the presented methods, a performance profile is used \cite{dolan2002}.

Performance profiles are designed to graphically compare both the speed and the robustness of algorithms across a test set. This is done by plotting, for each algorithm, the percentage of problems that are solved within a factor of the best solve time. For a more detailed description of performance profiles, see \cite{dolan2002}.

\begin{figure}[ht]
\centering
\subfigure[Maximum tolerance of 5\%]{
\includegraphics[height=2in]{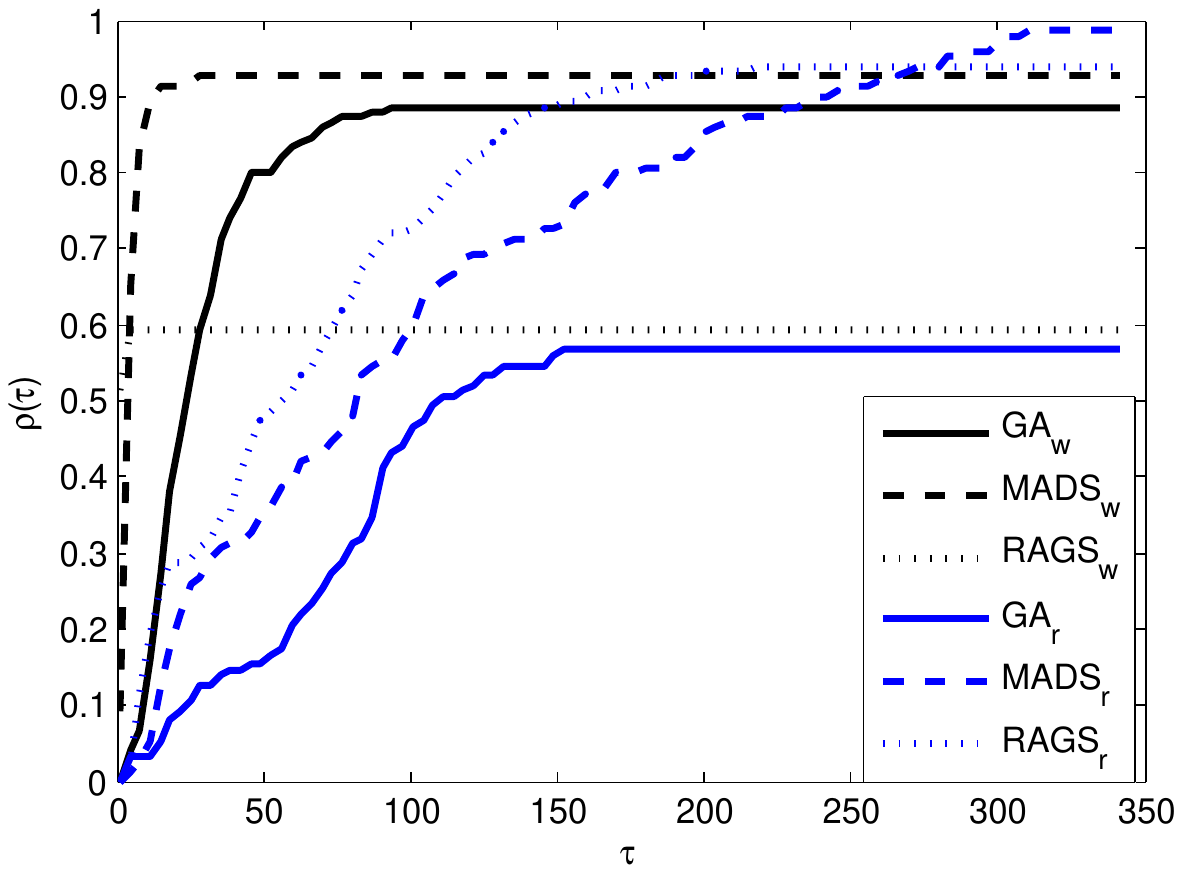}
\label{figure:profile5}
}
\subfigure[Maximum tolerance of 1\%]{
\includegraphics[height=2in]{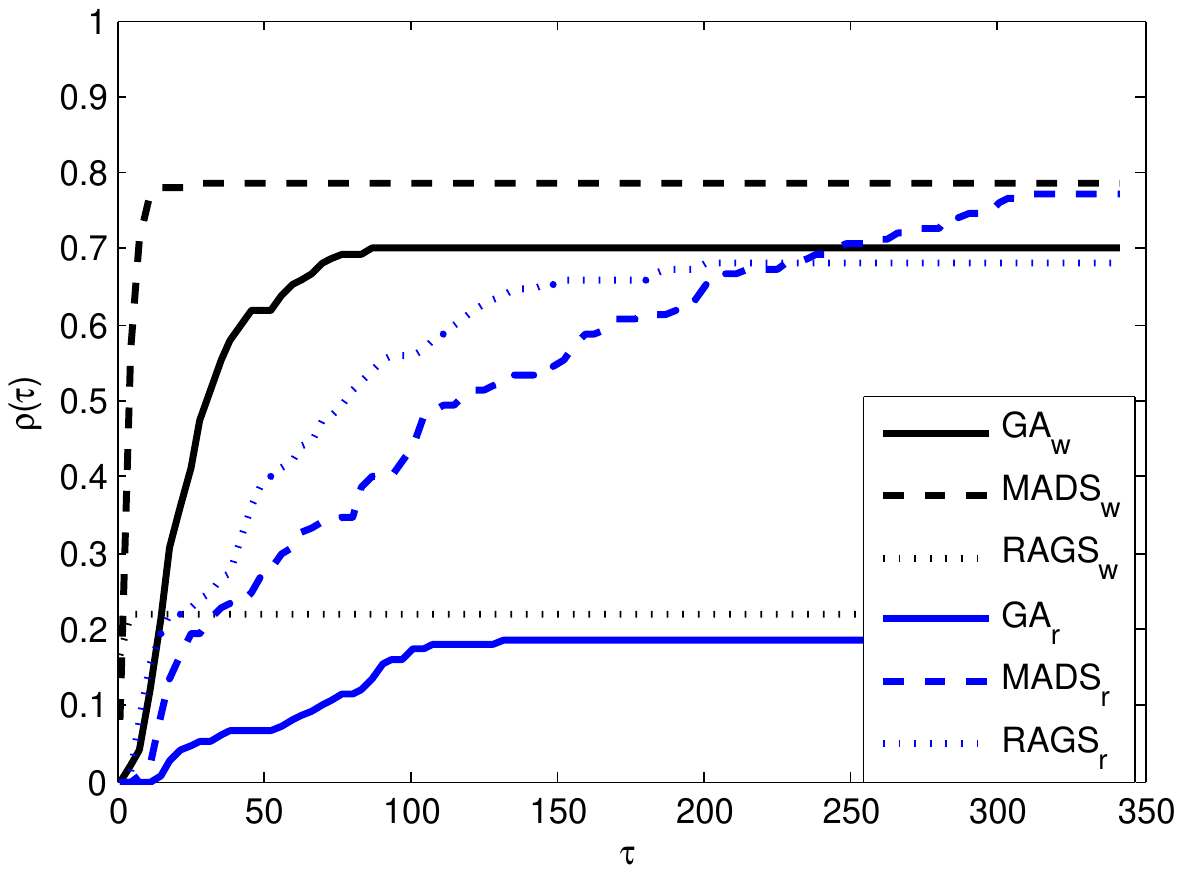}
\label{figure:profile1}
}
\caption{Performance profiles for GA$_w$, MADS$_w$, RAGS$_w$, GA$_r$, MADS$_r$, and RAGS$_r$.}\label{figure:profile}
\end{figure}

To calculate the performance profile, a definition of when a method ``solves" a specific problem is required. In this paper, a method is considered as a ``failed method" if the difference between the objective value obtained using the method in question and the best objective value obtained by any of the methods for that problem exceeds the defined allowable tolerance.
Performance profiles for the presented methods are plotted in Figure \ref{figure:profile} for allowable tolerances of 5\% and 1\%, respectively.

In Figure \ref{figure:profile5}, it is shown that for 5\% tolerance the maximum accuracy is obtained for MADS$_r$ and RAGS$_r$.  However, MADS$_w$ does extremely well, solving over 90\% of the problems, and takes only a fraction of the time portion to achieve this.  In Figure \ref{figure:profile1}, it is shown that for 1\% tolerance,  MADS$_w$ not only provides maximum accuracy, but also uses the least solving time.  Overall, it appears that the MADS algorithm using a warm-start procedure is well suited to solve these problems.

An interesting note occurs in comparing the warm-start with random initial points.  Warm-starting seems to give a small positive boost to RAGS, particularly in rate of convergence.  Conversely, RAGS$_r$ actually outperforms RAGS$_w$ in final accuracy.  This is likely because, unlike MADS and GA, RAGS has no embedded heuristics to break out of local minimizers.  This suggests that the warm-start locations, while good, are local minimizers of the next problem.  Finally, without warm-starting GA performs quite poorly.

\subsection{Number of Dampers}\label{number-of-dampers}

In this section, the effects of the number of dampers on the efficiency of the retrofitting system for one of the test problems are presented. As a multi-objective optimization study, fewer dampers and increased efficiency of the system are desired.  One of the key results in this section (and this research) is that, if a proper optimization is applied, then there is no need to place dampers on every story.  In fact, one can get optimal results by placing dampers on only a fraction of the total number of stories.

\begin{figure}[ht]
\centering
\subfigure[Material Set I]{
\includegraphics[height=2in]{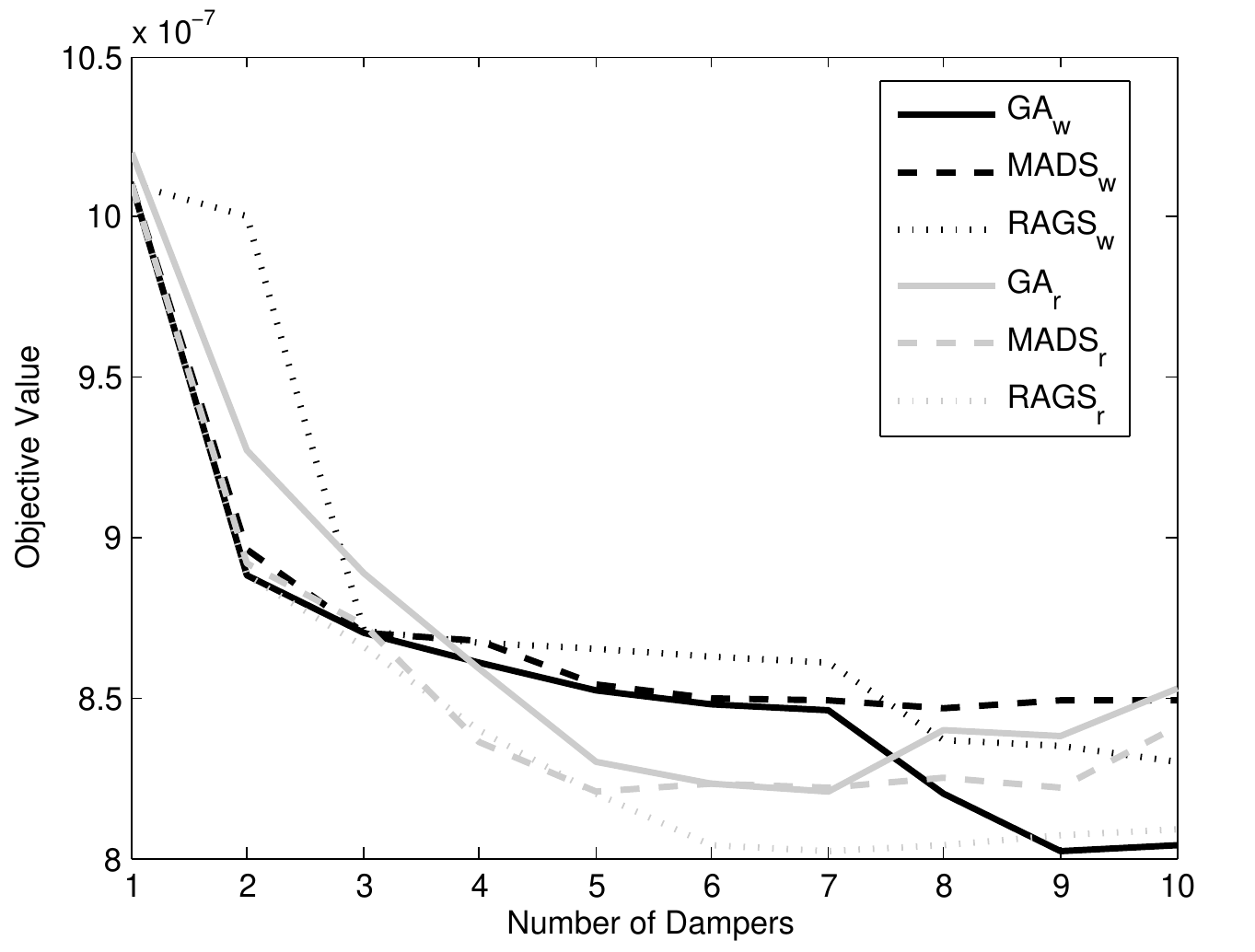}
\label{figure:dampersmat1}
}
\subfigure[Material Set II]{
\includegraphics[height=2in]{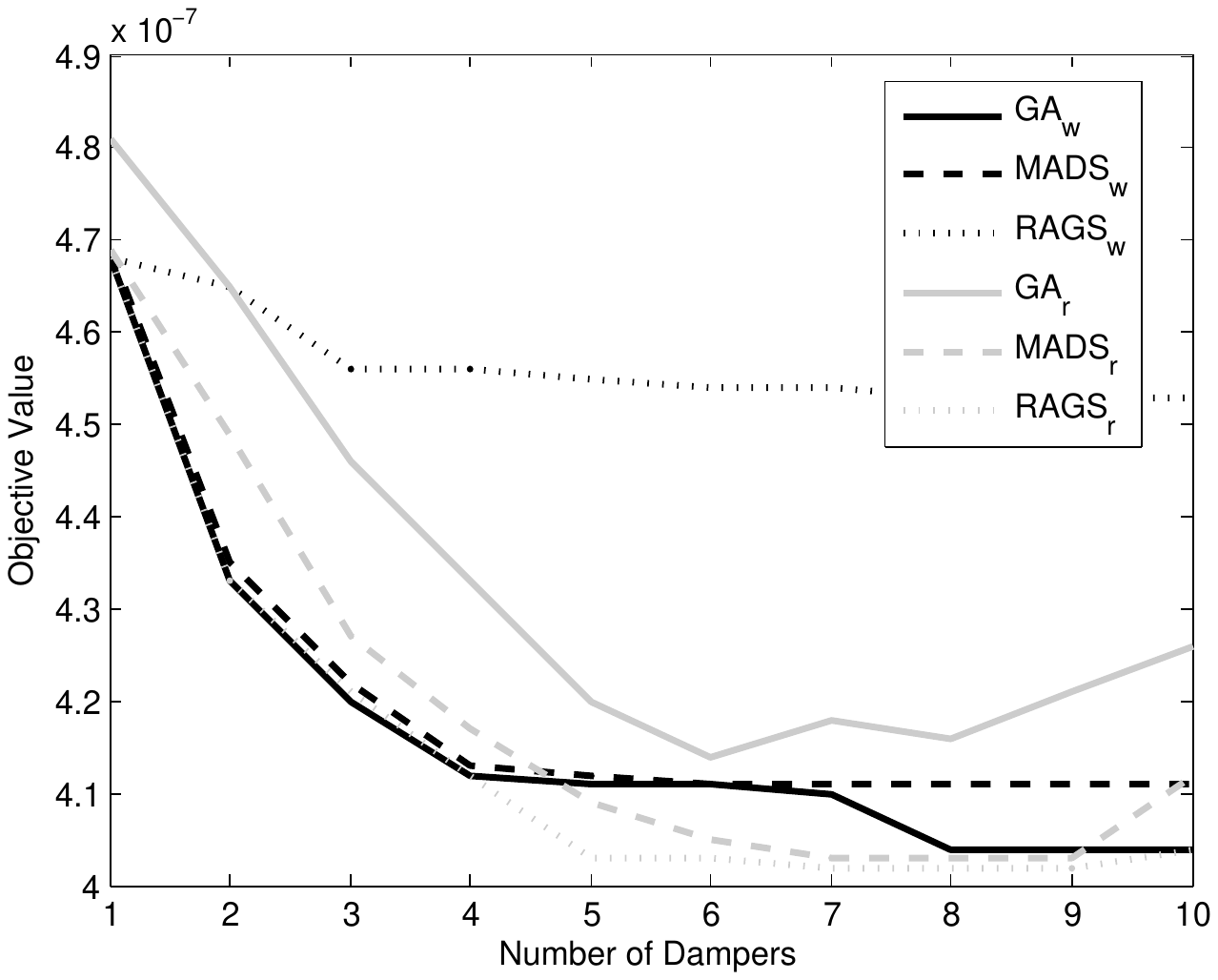}
\label{figure:dampersmat2}
}
\subfigure[Material Set III]{
\includegraphics[height=2in]{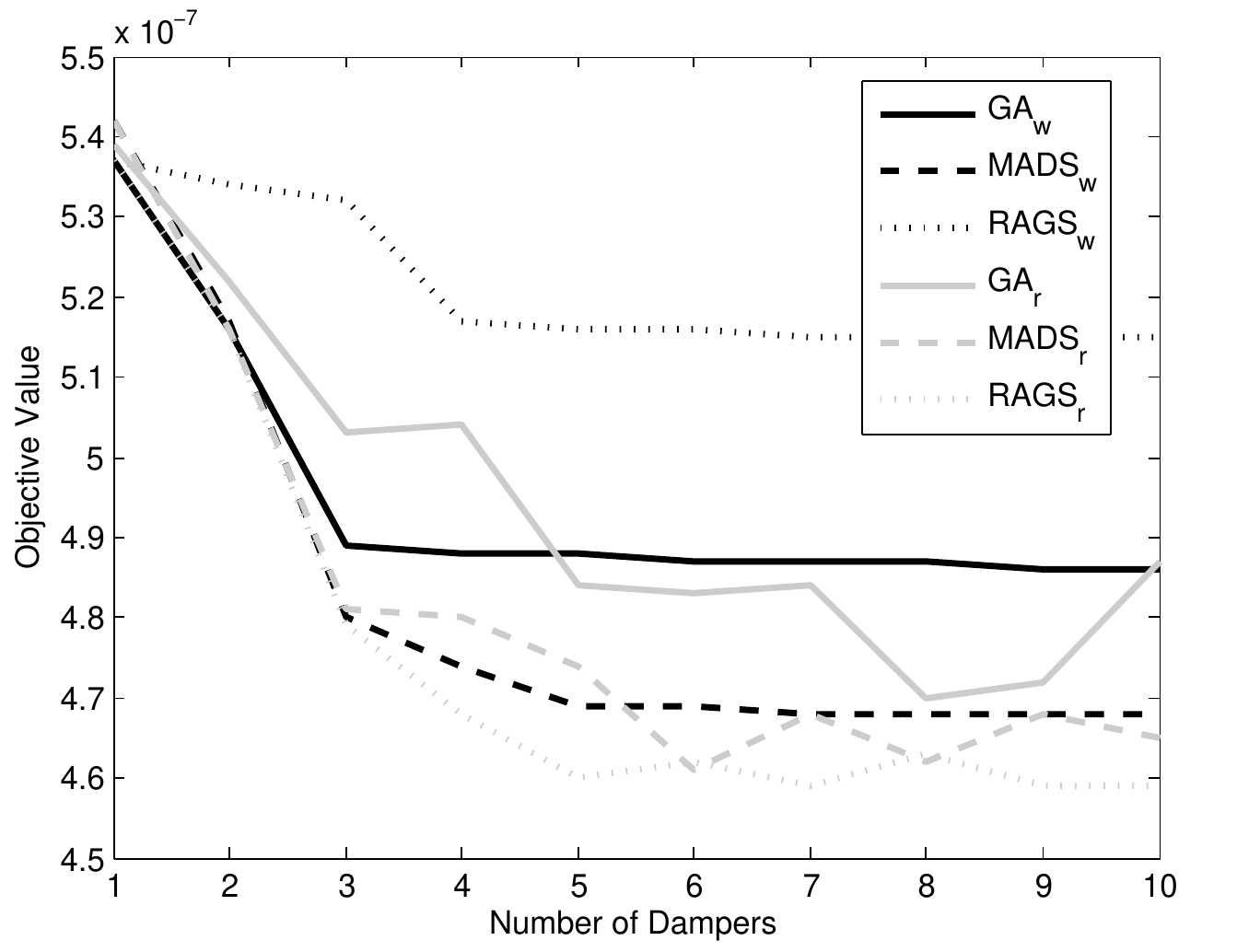}
\label{figure:dampersmat3}
}
\caption{Objective value for an increasing number of inserted dampers using Building Heights I.}
\label{figure:numdampers}
\end{figure}

To help visualize this result, Figure \ref{figure:numdampers} plots the number of dampers used against the optimal maximum inter-story drift achieved for Building Heights I and Material Sets I, II and III.  (Plots for other Building Heights look similar, and are available by contacting the corresponding author.)  As expected, the objective value generally decreases as the number of dampers are increased.  What is surprising is how rapidly the objective value decreases.  For Material Set I, Figure \ref{figure:dampersmat1}, optimal values are obtained using just 5 or 6 dampers for every algorithm except GA$_w$.  Similar trends occur in Figures \ref{figure:dampersmat2} and \ref{figure:dampersmat3}.

Another interesting note occurs in comparing the warm-start with random initial points.  As expected, the use of warm-start initial points means that the objective value never increases when the number of dampers is increased.  When random initial points are employed, this trend is not present, and indeed GA$_r$ does notably worse using 10 dampers than using just 5 or 6.  On the other hand, examining Figures \ref{figure:dampersmat2} and \ref{figure:dampersmat3}, it is seen that RAGS$_r$ stuck in a local minimizer that requires just 3 dampers.  So, warm-starting appears valuable, but only if the algorithm includes some heuristic to break free of local minimizers.

\begin{figure}[ht]
\centering
\subfigure[Exact Minimal Solution]{
\includegraphics[width=2in]{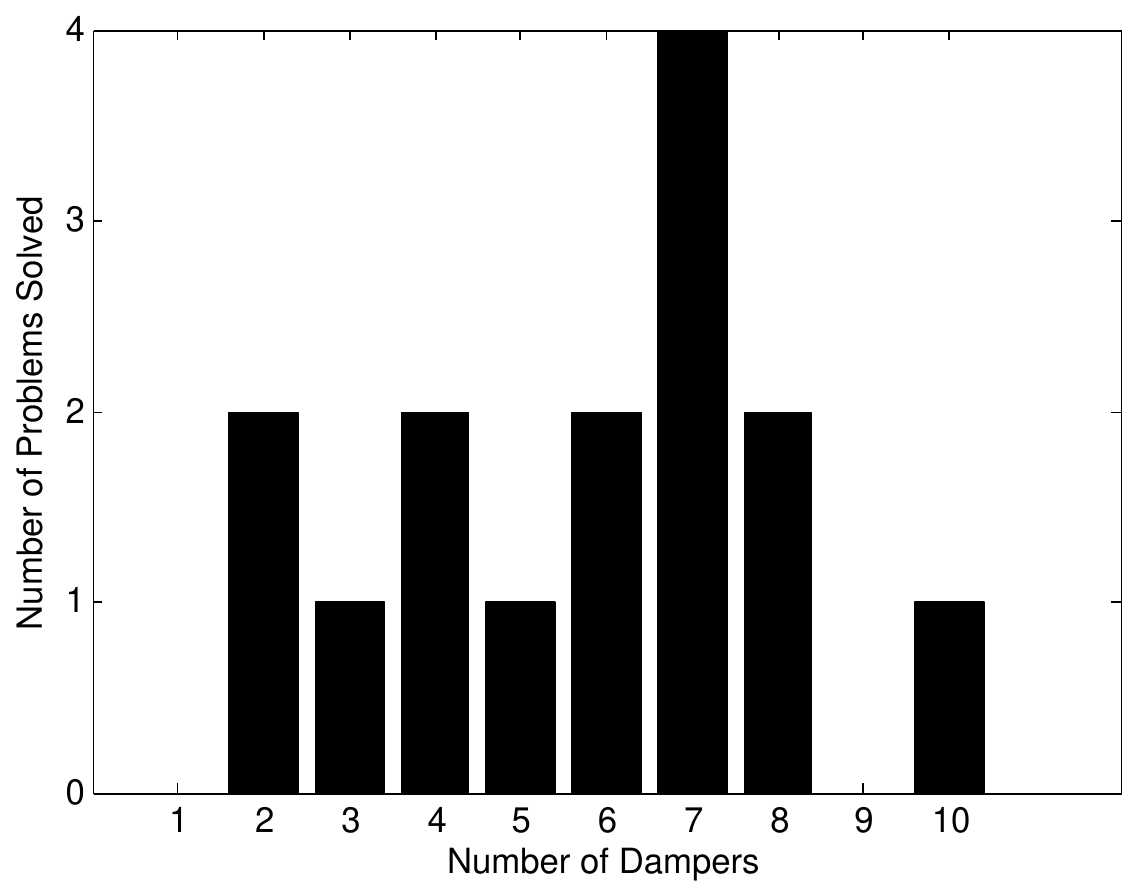}
\label{figure:damperhistogram0}
}
\subfigure[Tolerance of 1\%]{
\includegraphics[width=2in]{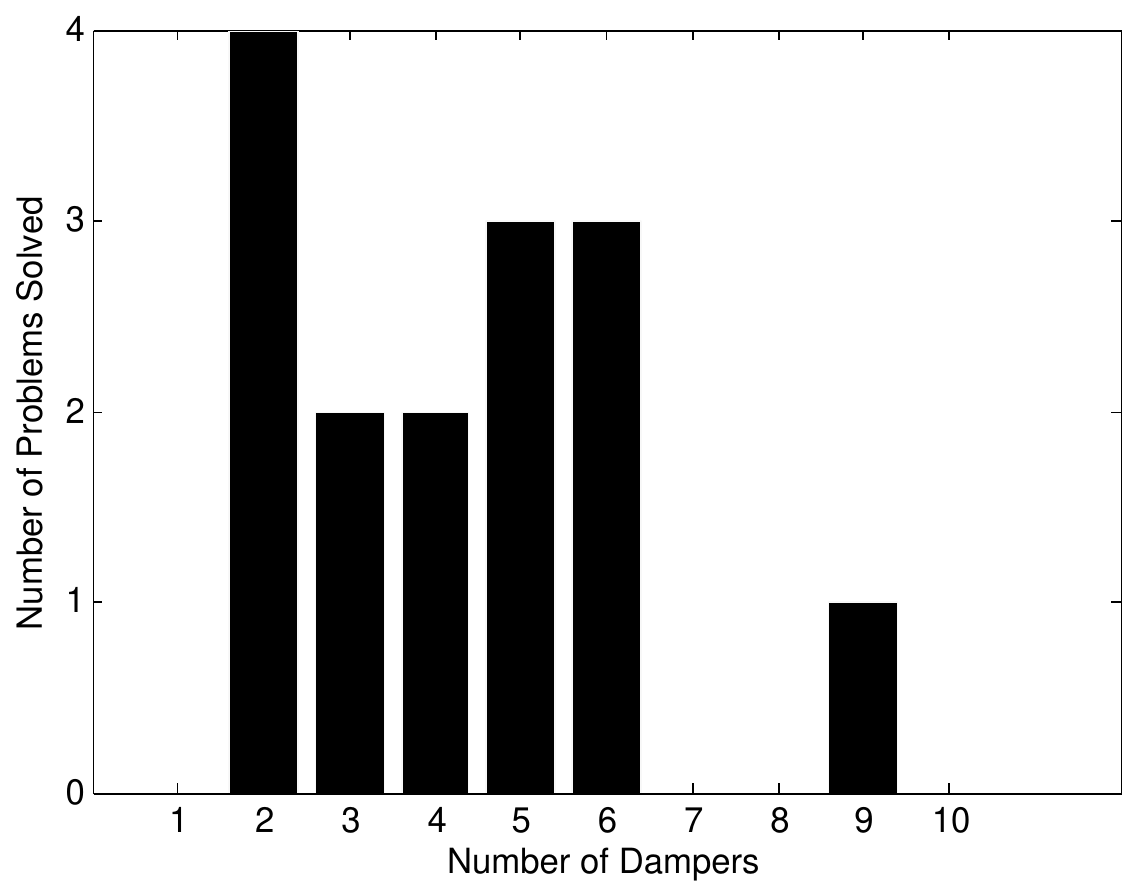}
\label{figure:damperhistogram1}
}
\caption{Number of dampers required to solve problems within a tolerance of the optimal solution.}\label{figure:damperhistogram}
\end{figure}

Figure \ref{figure:numdampers} inspires us to consider how many dampers are required by each of the algorithms to find an optimal solution for a single building.  For a fixed Mechanical Property and Building Height, the optimal solution is taken to be the overall lowest value found by all six algorithms given any number of dampers.  This yields 15 test problems (3 sets of Mechanical Properties and 5 sets of Building Heights).  The results are represented in Figure \ref{figure:damperhistogram}.

In Figure \ref{figure:damperhistogram0}, a histogram of the number of dampers used in the exact optimal solution for a fixed Mechanical Property and Building Height is provided.  While Figure \ref{figure:damperhistogram1}, provides a histogram of the minimum number of dampers used in order to minimize the objective within 1\% of the optimal solution. Examining Figure \ref{figure:damperhistogram0}, notice that only one problem requires 10 dampers to achieve the minimum value (this is Material Set III, Building Height 4).  More interestingly, in Figure \ref{figure:damperhistogram1} it is seen that the vast majority of problems are solved with less than 6 dampers, and many require just 2 or 3 dampers to achieve a solution within 1\% of the optimal solution.  It is worth noting that the problem that requires 9 dampers  to achieve a solution within 1\% of the optimal solution is also Material Set III, Building Height 4.

\subsection{Damper Configuration}

\begin{figure}[ht]
\centering
\subfigure[1 Damper]{
\includegraphics[width=2in]{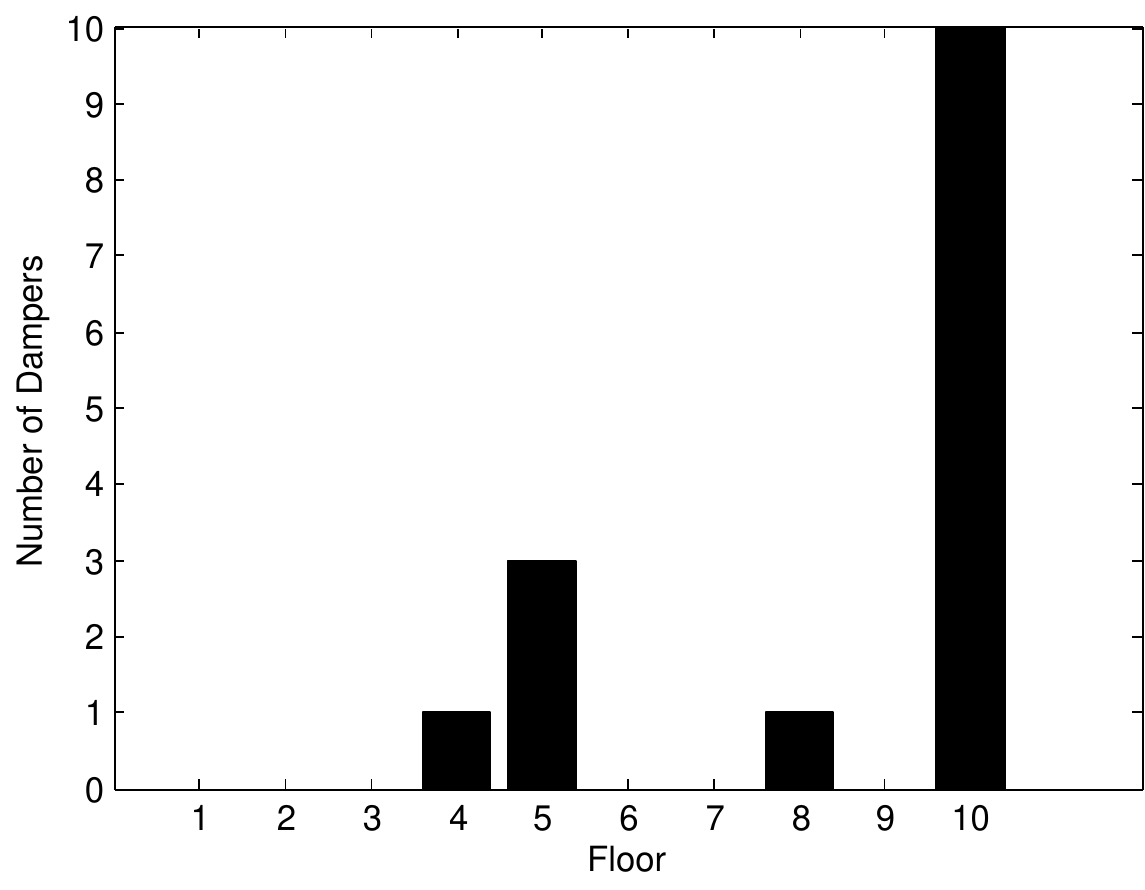}
\label{figure:damperlocations1}
}
\subfigure[2 Dampers]{
\includegraphics[width=2in]{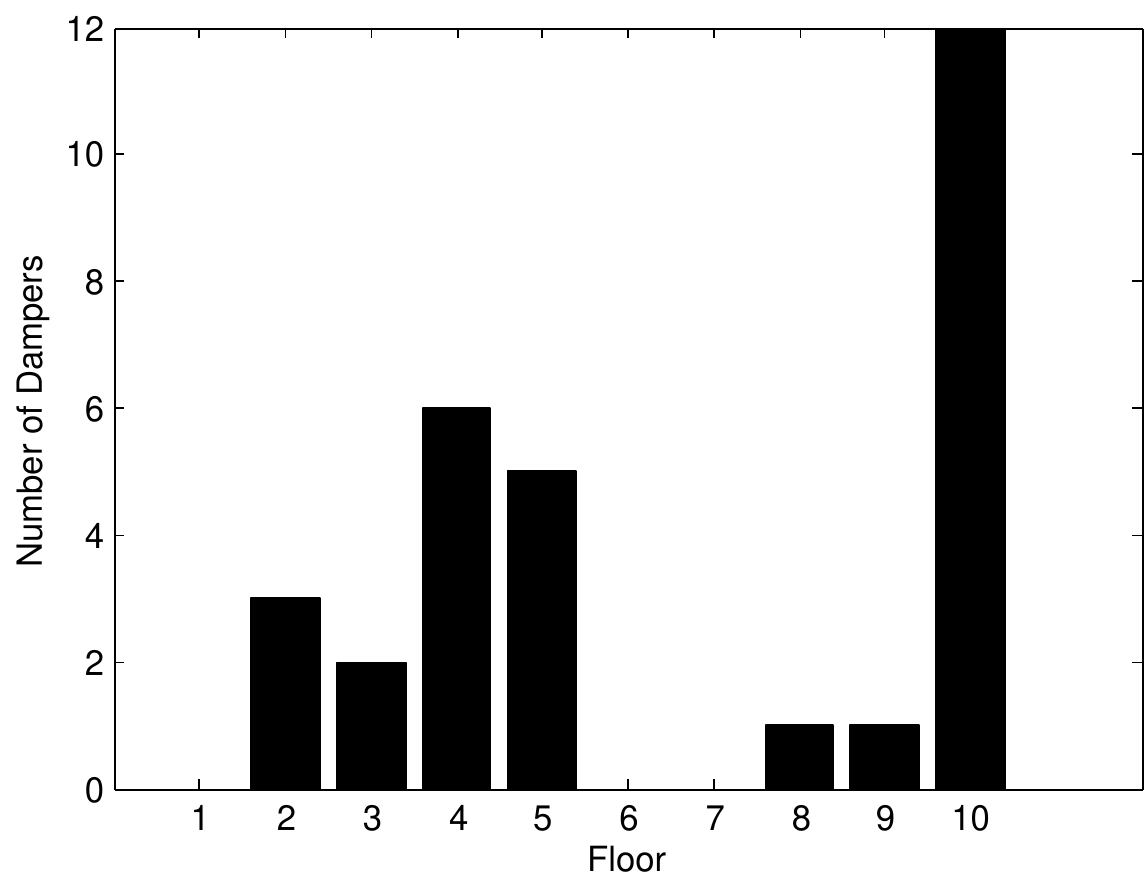}
\label{figure:damperlocations2}
}
\subfigure[5 Dampers]{
\includegraphics[width=2in]{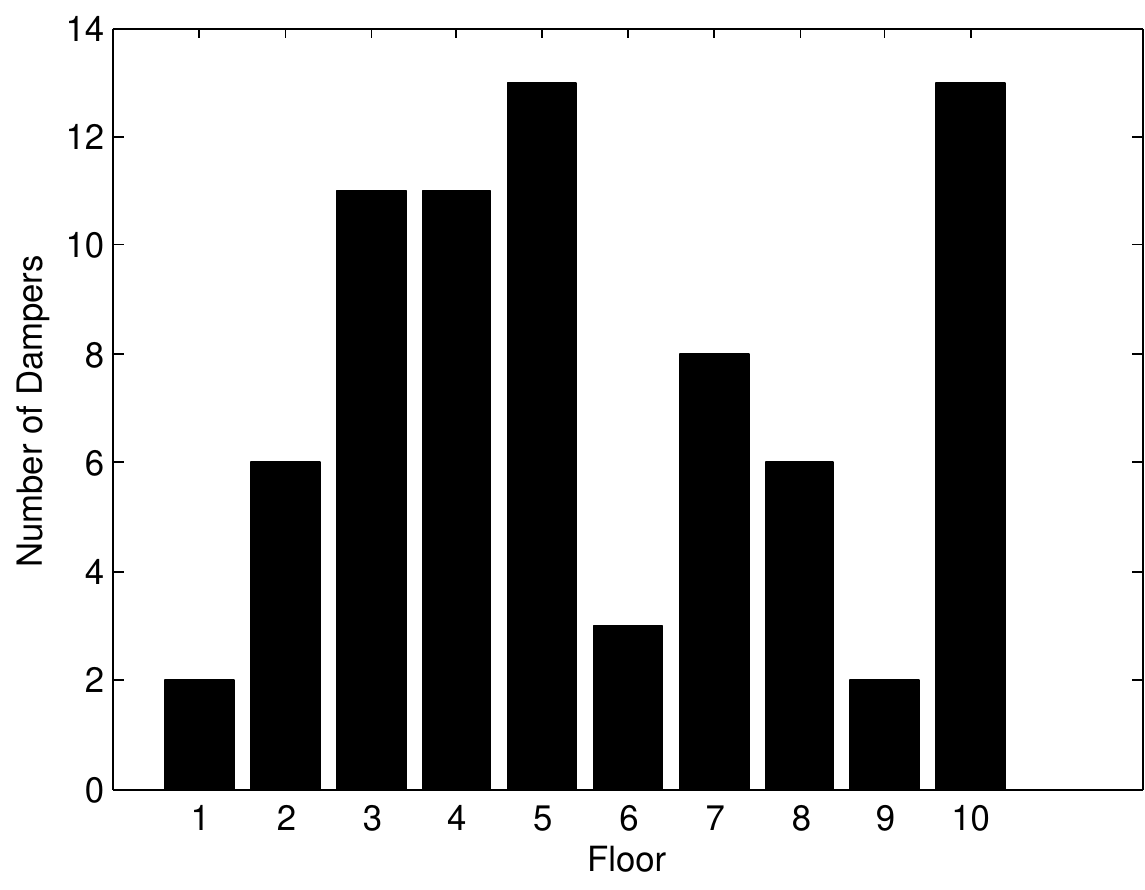}
\label{figure:damperlocations5}
}
\caption{Histograms of aggregate optimal damper locations for the 15 test buildings.}\label{figure:damperlocations}
\end{figure}

In Subsection \ref{number-of-dampers}, it is found that for most building combinations, the number of dampers required to solve problems within a tolerance of the optimal solution is less than half of the maximum number of dampers that could be inserted.  It is particularly interesting that several problems can be solved (with 1\% tolerance) using just 2 dampers.  A natural question at this point is, at which floors are dampers most commonly inserted?   Figure \ref{figure:damperlocations} examines this question, specifically looking at the cases when 1 damper, 2 dampers, and 5 dampers are inserted.

Figure \ref{figure:damperlocations} plots a histogram of the optimal damper locations for the 15 buildings.  Notice that if only one damper is used, then by far the most common location is on the tenth floor (as high as possible in the problem). Examining  Figures \ref{figure:damperlocations2} and \ref{figure:damperlocations5}, the pattern becomes less apparent.  The most popular location is always the tenth floor, but as more dampers are added, the locations become more scattered.  This emphasizes the importance of optimization and considering each building configuration uniquely.

\section{Conclusion}\label{Conclusions}

This paper presents a comprehensive optimization problem formulation and procedure that can be used to find the optimal configuration and mechanical properties of dampers for connected structures.
In particular, two adjacent buildings are considered as lumped mass models connected to each other using discrete viscous dampers. A pseudo excitation formula is used to generate an earthquake load in a frequency domain. Assuming a linear behaviour of the buildings (linear springs and linear viscous dampers), the dynamic response of the whole system is found. Using the dynamic response of the system, the desired objective function, i.e., the maximum inter-story drift, is calculated.

The optimization procedure consists of two parts including discrete and continuous optimizations. An outer-loop (discrete optimization algorithm) finds the best configuration of a limited number of dampers between two buildings; an inner-loop (continuous optimization algorithm) finds the optimal damper coefficients of the dampers. Three different algorithms (GA, MADS and RAGS) for the continuous optimization problem are considered, each using a random initial point and a warm-start initial point.  In order to compare speed and robustness of these non-gradient based methods, $150$ test problems were generated and solved via these three methods. Results showed that MADS using a warm-start initial point is quite fast and robust.

Furthermore, the efficiency of the retrofitting system with respect to the number of dampers used was investigated. In \cite{opt_loc_bigdeli2012}, it is shown that when assuming equal damper coefficients, increasing the number of dampers may exacerbate the dynamic behaviour of the buildings. When the assumption of equal damper coefficients is removed, it was observed that although increasing the number of dampers no longer exacerbates the dynamic behaviour of the system, there is nonetheless a threshold after which increasing the number of dampers provides little benefit to the system.  Using 15 test problems, it was found that in most cases, the optimal behaviour of a seismic retrofit can be achieved within 1\% using 4 or less dampers, and only one problem required more than 6 dampers.  This represents a significant saving in material and overall cost of retrofitting.

Finally, it is worth mentioning that a very similar bi-level optimization procedure as presented can be followed for different types of damper connectors, such as MR dampers, friction dampers and so on. In these cases, the only element that changes is the simulation core of the problem; the same discrete and continuous optimization algorithms can be used.  Furthermore, it should be clear that, while this paper focused on minimizing the maximum inter-story drift, the techniques within this paper can easily be adapted to any objective function.

\bibliographystyle{plain}	
\bibliography{Refs}

\newpage

\section{Appendix A}\label{data-tables}

\begin{table}[h!]\scriptsize
  \centering
  \caption{Number of function calls and objective values found for Material Set I for algorithms GA$_r$, MADS$_r$, and RAGS$_r$.}
    \begin{tabular}{@{}cccccccccc@{}}
    \toprule
     && \multicolumn{3}{c}{Number of function calls}  && \multicolumn{3}{c}{Objective value}  \\[-0.2mm]
    \midrule
    Case  & nd    & GA$_r$    & MADS$_r$  & RAGS$_r$  &       & GA$_r$    & MADS$_r$  & RAGS$_r$ \\[-0.2mm]
    \midrule
   1	&	1	&	380	&	272	&	203	&	&	1.02E-06	&	1.01E-06	&	1.01E-06	\\
	&	2	&	1030	&	930	&	837	&	&	9.27E-07	&	8.92E-07	&	8.88E-07	\\
	&	3	&	1870	&	1849	&	1455	&	&	8.89E-07	&	8.73E-07	&	8.66E-07	\\
	&	4	&	2850	&	3014	&	2297	&	&	8.59E-07	&	8.36E-07	&	8.40E-07	\\
	&	5	&	3900	&	4445	&	3056	&	&	8.30E-07	&	8.21E-07	&	8.20E-07	\\
	&	6	&	4950	&	6224	&	3892	&	&	8.23E-07	&	8.23E-07	&	8.04E-07	\\
	&	7	&	5930	&	8192	&	4746	&	&	8.21E-07	&	8.22E-07	&	8.02E-07	\\
	&	8	&	6770	&	10282	&	5600	&	&	8.40E-07	&	8.25E-07	&	8.04E-07	\\
	&	9	&	7400	&	11327	&	6212	&	&	8.38E-07	&	8.22E-07	&	8.07E-07	\\
	&	10	&	350	&	695	&	552	&	&	8.53E-07	&	8.41E-07	&	8.09E-07	\\
\midrule																
2	&	1	&	380	&	279	&	170	&	&	4.36E-06	&	4.33E-06	&	4.33E-06	\\
	&	2	&	1100	&	873	&	571	&	&	1.79E-06	&	1.72E-06	&	1.72E-06	\\
	&	3	&	2000	&	1883	&	1661	&	&	1.90E-06	&	1.72E-06	&	1.72E-06	\\
	&	4	&	2980	&	3680	&	2965	&	&	1.87E-06	&	1.72E-06	&	1.71E-06	\\
	&	5	&	4055	&	5995	&	4360	&	&	1.76E-06	&	1.72E-06	&	1.73E-06	\\
	&	6	&	5105	&	8999	&	6376	&	&	1.77E-06	&	1.72E-06	&	1.73E-06	\\
	&	7	&	6155	&	13127	&	9192	&	&	1.84E-06	&	1.71E-06	&	1.73E-06	\\
	&	8	&	6995	&	16589	&	11644	&	&	1.88E-06	&	1.71E-06	&	1.74E-06	\\
	&	9	&	7625	&	18710	&	13511	&	&	1.89E-06	&	1.73E-06	&	1.74E-06	\\
	&	10	&	350	&	1511	&	1208	&	&	1.96E-06	&	1.71E-06	&	1.77E-06	\\
\midrule																
3	&	1	&	370	&	266	&	164	&	&	1.67E-06	&	1.59E-06	&	1.59E-06	\\
	&	2	&	1020	&	729	&	724	&	&	1.62E-06	&	1.53E-06	&	1.52E-06	\\
	&	3	&	1860	&	1367	&	1719	&	&	1.56E-06	&	1.52E-06	&	1.51E-06	\\
	&	4	&	2840	&	2030	&	2814	&	&	1.54E-06	&	1.51E-06	&	1.51E-06	\\
	&	5	&	3890	&	2933	&	3911	&	&	1.54E-06	&	1.52E-06	&	1.52E-06	\\
	&	6	&	4940	&	4133	&	5228	&	&	1.55E-06	&	1.52E-06	&	1.52E-06	\\
	&	7	&	5920	&	5780	&	6443	&	&	1.54E-06	&	1.52E-06	&	1.52E-06	\\
	&	8	&	6760	&	8014	&	7794	&	&	1.54E-06	&	1.52E-06	&	1.53E-06	\\
	&	9	&	7390	&	9667	&	8845	&	&	1.54E-06	&	1.54E-06	&	1.53E-06	\\
	&	10	&	350	&	1812	&	513	&	&	1.57E-06	&	1.52E-06	&	1.54E-06	\\
\midrule																
4	&	1	&	405	&	322	&	154	&	&	8.24E-06	&	8.21E-06	&	8.19E-06	\\
	&	2	&	1065	&	1045	&	640	&	&	5.87E-06	&	5.58E-06	&	5.56E-06	\\
	&	3	&	1920	&	2190	&	1247	&	&	5.57E-06	&	5.48E-06	&	5.47E-06	\\
	&	4	&	2940	&	3962	&	1846	&	&	5.57E-06	&	5.45E-06	&	5.44E-06	\\
	&	5	&	3990	&	6221	&	2598	&	&	5.47E-06	&	5.44E-06	&	5.43E-06	\\
	&	6	&	5040	&	9050	&	3402	&	&	5.46E-06	&	5.43E-06	&	5.42E-06	\\
	&	7	&	6020	&	11796	&	4167	&	&	5.46E-06	&	5.42E-06	&	5.42E-06	\\
	&	8	&	6860	&	14535	&	4941	&	&	5.44E-06	&	5.42E-06	&	5.43E-06	\\
	&	9	&	7490	&	17575	&	5666	&	&	5.53E-06	&	5.42E-06	&	5.45E-06	\\
	&	10	&	350	&	1933	&	487	&	&	5.66E-06	&	5.43E-06	&	5.46E-06	\\
\midrule																
5	&	1	&	385	&	312	&	153	&	&	5.96E-06	&	5.94E-06	&	5.94E-06	\\
	&	2	&	1045	&	1010	&	756	&	&	5.86E-06	&	5.86E-06	&	5.87E-06	\\
	&	3	&	1930	&	2097	&	1766	&	&	5.94E-06	&	5.86E-06	&	5.87E-06	\\
	&	4	&	2930	&	4345	&	3424	&	&	5.89E-06	&	5.86E-06	&	5.87E-06	\\
	&	5	&	3980	&	7131	&	5166	&	&	6.17E-06	&	5.86E-06	&	5.88E-06	\\
	&	6	&	5030	&	11821	&	7045	&	&	5.96E-06	&	5.86E-06	&	5.89E-06	\\
	&	7	&	6010	&	16077	&	8844	&	&	6.22E-06	&	5.86E-06	&	5.89E-06	\\
	&	8	&	6850	&	21024	&	10533	&	&	5.97E-06	&	5.86E-06	&	5.91E-06	\\
	&	9	&	7480	&	25831	&	11646	&	&	6.14E-06	&	5.86E-06	&	5.92E-06	\\
	&	10	&	350	&	3545	&	820	&	&	6.40E-06	&	5.86E-06	&	5.92E-06	\\
    \bottomrule
    \end{tabular}%
  \label{table:dataMatIrandom}%
\end{table}%

\begin{table}[h!]\scriptsize
  \centering
  \caption{Number of function calls and objective values found for Material Set I for algorithms GA$_w$, MADS$_w$, and RAGS$_w$.}
    \begin{tabular}{@{}cccccccccc@{}}
    \toprule
     && \multicolumn{3}{c}{Number of function calls}  && \multicolumn{3}{c}{Objective value}  \\[-0.2mm]
    \midrule    Case  & nd    & GA$_w$    & MADS$_w$  & RAGS$_w$  &       & GA$_w$    & MADS$_w$  & RAGS$_w$ \\[-0.2mm]
    \midrule
1	&	1	&	1040	&	35	&	19	&&	1.01E-06	&	1.01E-06	&	1.01E-06	\\
	&	2	&	1060	&	68	&	32	&&	8.88E-07	&	8.96E-07	&	1.00E-06	\\
	&	3	&	1040	&	110	&	242	&&	8.70E-07	&	8.70E-07	&	8.71E-07	\\
	&	4	&	1060	&	48	&	71	&&	8.61E-07	&	8.68E-07	&	8.67E-07	\\
	&	5	&	1040	&	326	&	44	&&	8.52E-07	&	8.54E-07	&	8.65E-07	\\
	&	6	&	1040	&	534	&	84	&&	8.48E-07	&	8.50E-07	&	8.63E-07	\\
	&	7	&	1040	&	160	&	113	&&	8.46E-07	&	8.49E-07	&	8.61E-07	\\
	&	8	&	1040	&	524	&	686	&&	8.20E-07	&	8.47E-07	&	8.37E-07	\\
	&	9	&	1040	&	205	&	144	&&	8.02E-07	&	8.49E-07	&	8.35E-07	\\
	&	10	&	1040	&	92	&	281	&&	8.04E-07	&	8.49E-07	&	8.30E-07	\\
\midrule
2	&	1	&	1040	&	26	&	23	&&	4.33E-06	&	4.33E-06	&	4.33E-06	\\
	&	2	&	1040	&	76	&	151	&&	1.72E-06	&	1.73E-06	&	1.72E-06	\\
	&	3	&	1040	&	35	&	51	&&	1.72E-06	&	1.72E-06	&	1.72E-06	\\
	&	4	&	1040	&	95	&	45	&&	1.72E-06	&	1.72E-06	&	1.72E-06	\\
	&	5	&	1040	&	317	&	40	&&	1.72E-06	&	1.72E-06	&	1.72E-06	\\
	&	6	&	1060	&	182	&	57	&&	1.72E-06	&	1.72E-06	&	1.72E-06	\\
	&	7	&	1040	&	254	&	78	&&	1.72E-06	&	1.72E-06	&	1.73E-06	\\
	&	8	&	1040	&	1188	&	143	&&	1.72E-06	&	1.72E-06	&	1.74E-06	\\
	&	9	&	1040	&	391	&	383	&&	1.72E-06	&	1.72E-06	&	1.74E-06	\\
	&	10	&	1040	&	351	&	214	&&	1.72E-06	&	1.72E-06	&	1.74E-06	\\
\midrule
3	&	1	&	1040	&	31	&	31	&&	1.59E-06	&	1.59E-06	&	1.59E-06	\\
	&	2	&	1040	&	60	&	16	&&	1.53E-06	&	1.53E-06	&	1.58E-06	\\
	&	3	&	1040	&	85	&	32	&&	1.52E-06	&	1.52E-06	&	1.58E-06	\\
	&	4	&	1040	&	107	&	41	&&	1.52E-06	&	1.52E-06	&	1.57E-06	\\
	&	5	&	1040	&	256	&	50	&&	1.52E-06	&	1.52E-06	&	1.57E-06	\\
	&	6	&	1040	&	133	&	76	&&	1.52E-06	&	1.52E-06	&	1.57E-06	\\
	&	7	&	1040	&	204	&	56	&&	1.51E-06	&	1.52E-06	&	1.56E-06	\\
	&	8	&	1040	&	179	&	97	&&	1.51E-06	&	1.52E-06	&	1.56E-06	\\
	&	9	&	1040	&	447	&	123	&&	1.51E-06	&	1.52E-06	&	1.56E-06	\\
	&	10	&	1041	&	261	&	127	&&	1.51E-06	&	1.52E-06	&	1.56E-06	\\
\midrule
4	&	1	&	1041	&	32	&	19	&&	8.19E-06	&	8.21E-06	&	8.19E-06	\\
	&	2	&	1041	&	66	&	70	&&	5.56E-06	&	5.60E-06	&	5.56E-06	\\
	&	3	&	1041	&	123	&	33	&&	5.47E-06	&	5.49E-06	&	5.55E-06	\\
	&	4	&	1041	&	168	&	41	&&	5.45E-06	&	5.45E-06	&	5.54E-06	\\
	&	5	&	1041	&	242	&	62	&&	5.43E-06	&	5.43E-06	&	5.54E-06	\\
	&	6	&	1041	&	176	&	39	&&	5.42E-06	&	5.42E-06	&	5.53E-06	\\
	&	7	&	1041	&	133	&	60	&&	5.42E-06	&	5.42E-06	&	5.52E-06	\\
	&	8	&	1041	&	595	&	64	&&	5.41E-06	&	5.42E-06	&	5.52E-06	\\
	&	9	&	1041	&	133	&	51	&&	5.42E-06	&	5.42E-06	&	5.53E-06	\\
	&	10	&	1041	&	184	&	106	&&	5.42E-06	&	5.42E-06	&	5.54E-06	\\
\midrule
5	&	1	&	1041	&	31	&	23	&&	5.94E-06	&	5.94E-06	&	5.94E-06	\\
	&	2	&	1041	&	54	&	28	&&	5.86E-06	&	5.86E-06	&	5.89E-06	\\
	&	3	&	1041	&	252	&	44	&&	5.86E-06	&	5.86E-06	&	5.87E-06	\\
	&	4	&	1041	&	301	&	42	&&	5.86E-06	&	5.86E-06	&	5.87E-06	\\
	&	5	&	1041	&	202	&	43	&&	5.86E-06	&	5.86E-06	&	5.87E-06	\\
	&	6	&	1041	&	196	&	39	&&	5.86E-06	&	5.86E-06	&	5.87E-06	\\
	&	7	&	1041	&	59	&	76	&&	5.86E-06	&	5.86E-06	&	5.87E-06	\\
	&	8	&	1041	&	98	&	73	&&	5.86E-06	&	5.86E-06	&	5.88E-06	\\
	&	9	&	1041	&	789	&	84	&&	5.86E-06	&	5.86E-06	&	5.90E-06	\\
	&	10	&	1041	&	1153	&	113	&&	5.86E-06	&	5.86E-06	&	5.91E-06	\\
    \bottomrule
    \end{tabular}%
  \label{table:dataMatIwarm}%
\end{table}%

\begin{table}[htbp]\scriptsize
  \centering
  \caption{Number of function calls and objective values found for Material Set II for algorithms GA$_r$, MADS$_r$, and RAGS$_r$.}
    \begin{tabular}{@{}cccccccccc@{}}
    \toprule
     && \multicolumn{3}{c}{Number of function calls}  && \multicolumn{3}{c}{Objective value}  \\[-0.2mm]
    \midrule
    Case  & nd    & GA    & MADS  & RAGS  &       & GA    & MADS  & RAGS \\[-0.2mm]
    \midrule
    1	&	1	&	425	&	298	&	153	&	&	4.81E-07	&	4.69E-07	&	4.68E-07	\\
	&	2	&	1055	&	871	&	575	&	&	4.65E-07	&	4.49E-07	&	4.33E-07	\\
	&	3	&	1895	&	1884	&	1421	&	&	4.46E-07	&	4.27E-07	&	4.21E-07	\\
	&	4	&	2875	&	3406	&	2264	&	&	4.33E-07	&	4.17E-07	&	4.12E-07	\\
	&	5	&	3925	&	5119	&	3236	&	&	4.20E-07	&	4.09E-07	&	4.03E-07	\\
	&	6	&	4975	&	6565	&	4496	&	&	4.14E-07	&	4.05E-07	&	4.03E-07	\\
	&	7	&	5955	&	8001	&	5680	&	&	4.18E-07	&	4.03E-07	&	4.02E-07	\\
	&	8	&	6795	&	10055	&	6503	&	&	4.16E-07	&	4.03E-07	&	4.02E-07	\\
	&	9	&	7425	&	11942	&	7259	&	&	4.21E-07	&	4.03E-07	&	4.02E-07	\\
	&	10	&	350	&	930	&	460	&	&	4.26E-07	&	4.12E-07	&	4.04E-07	\\
\midrule																
2	&	1	&	355	&	324	&	196	&	&	2.43E-07	&	2.43E-07	&	2.43E-07	\\
	&	2	&	1025	&	1053	&	550	&	&	2.41E-07	&	2.37E-07	&	2.36E-07	\\
	&	3	&	1895	&	2203	&	1036	&	&	2.38E-07	&	2.37E-07	&	2.36E-07	\\
	&	4	&	2875	&	3460	&	1429	&	&	2.37E-07	&	2.36E-07	&	2.36E-07	\\
	&	5	&	3925	&	4481	&	1894	&	&	2.37E-07	&	2.37E-07	&	2.37E-07	\\
	&	6	&	4975	&	5281	&	2263	&	&	2.37E-07	&	2.38E-07	&	2.37E-07	\\
	&	7	&	5955	&	6159	&	2810	&	&	2.37E-07	&	2.36E-07	&	2.37E-07	\\
	&	8	&	6795	&	7178	&	3288	&	&	2.37E-07	&	2.38E-07	&	2.38E-07	\\
	&	9	&	7425	&	7632	&	3780	&	&	2.37E-07	&	2.40E-07	&	2.38E-07	\\
	&	10	&	350	&	172	&	307	&	&	2.42E-07	&	2.41E-07	&	2.38E-07	\\
\midrule																
3	&	1	&	410	&	308	&	123	&	&	5.35E-07	&	4.99E-07	&	4.99E-07	\\
	&	2	&	1060	&	1036	&	611	&	&	4.68E-07	&	3.32E-07	&	3.22E-07	\\
	&	3	&	1945	&	2130	&	1740	&	&	3.82E-07	&	3.24E-07	&	3.19E-07	\\
	&	4	&	2945	&	3689	&	3284	&	&	3.50E-07	&	3.21E-07	&	3.19E-07	\\
	&	5	&	4020	&	5121	&	4374	&	&	3.42E-07	&	3.20E-07	&	3.18E-07	\\
	&	6	&	5070	&	7196	&	5776	&	&	3.61E-07	&	3.13E-07	&	3.17E-07	\\
	&	7	&	6050	&	9158	&	7933	&	&	3.33E-07	&	3.11E-07	&	3.21E-07	\\
	&	8	&	6930	&	10932	&	10490	&	&	3.35E-07	&	3.15E-07	&	3.16E-07	\\
	&	9	&	7560	&	12724	&	12279	&	&	3.41E-07	&	3.20E-07	&	3.15E-07	\\
	&	10	&	350	&	1401	&	1196	&	&	3.54E-07	&	3.24E-07	&	3.13E-07	\\
\midrule																
4	&	1	&	390	&	299	&	409	&	&	3.89E-06	&	3.86E-06	&	3.86E-06	\\
	&	2	&	1060	&	954	&	851	&	&	3.50E-06	&	2.74E-06	&	2.74E-06	\\
	&	3	&	1945	&	1889	&	1663	&	&	2.74E-06	&	2.68E-06	&	2.72E-06	\\
	&	4	&	2945	&	3383	&	2903	&	&	2.75E-06	&	2.66E-06	&	2.70E-06	\\
	&	5	&	3995	&	5003	&	4040	&	&	2.77E-06	&	2.64E-06	&	2.71E-06	\\
	&	6	&	5075	&	7171	&	5382	&	&	2.73E-06	&	2.62E-06	&	2.70E-06	\\
	&	7	&	6055	&	8987	&	6510	&	&	2.75E-06	&	2.62E-06	&	2.71E-06	\\
	&	8	&	6895	&	10931	&	7559	&	&	2.79E-06	&	2.61E-06	&	2.70E-06	\\
	&	9	&	7525	&	12645	&	8550	&	&	2.77E-06	&	2.61E-06	&	2.71E-06	\\
	&	10	&	350	&	981	&	573	&	&	2.78E-06	&	2.61E-06	&	2.71E-06	\\
\midrule																
5	&	1	&	415	&	299	&	199	&	&	2.06E-06	&	1.93E-06	&	1.93E-06	\\
	&	2	&	1085	&	945	&	699	&	&	1.65E-06	&	1.24E-06	&	1.24E-06	\\
	&	3	&	1940	&	2047	&	1424	&	&	1.52E-06	&	1.23E-06	&	1.23E-06	\\
	&	4	&	2940	&	3455	&	2284	&	&	1.31E-06	&	1.22E-06	&	1.22E-06	\\
	&	5	&	3990	&	6590	&	3236	&	&	1.22E-06	&	1.21E-06	&	1.22E-06	\\
	&	6	&	5040	&	9516	&	4253	&	&	1.25E-06	&	1.21E-06	&	1.22E-06	\\
	&	7	&	6020	&	12164	&	5278	&	&	1.26E-06	&	1.21E-06	&	1.22E-06	\\
	&	8	&	6860	&	15396	&	6126	&	&	1.25E-06	&	1.21E-06	&	1.22E-06	\\
	&	9	&	7490	&	18878	&	6842	&	&	1.27E-06	&	1.21E-06	&	1.22E-06	\\
	&	10	&	350	&	1422	&	421	&	&	1.26E-06	&	1.21E-06	&	1.23E-06	\\
    \bottomrule
    \end{tabular}%
  \label{table:dataMatIIrandom}%
\end{table}%

\begin{table}[h!]\scriptsize
  \centering
  \caption{Number of function calls and objective values found for Material Set II for algorithms GA$_w$, MADS$_w$, and RAGS$_w$.}
    \begin{tabular}{@{}cccccccccc@{}}
    \toprule
     && \multicolumn{3}{c}{Number of function calls}  && \multicolumn{3}{c}{Objective value}  \\[-0.2mm]
    \midrule    Case  & nd    & GA$_w$    & MADS$_w$  & RAGS$_w$  &       & GA$_w$    & MADS$_w$  & RAGS$_w$ \\[-0.2mm]
    \midrule
1	&	1	&	1041	&	25	&	18	&&	4.68E-07	&	4.69E-07	&	4.68E-07	\\
	&	2	&	1421	&	55	&	34	&&	4.33E-07	&	4.35E-07	&	4.65E-07	\\
	&	3	&	1041	&	85	&	320	&&	4.20E-07	&	4.22E-07	&	4.56E-07	\\
	&	4	&	1141	&	222	&	41	&&	4.12E-07	&	4.13E-07	&	4.56E-07	\\
	&	5	&	1041	&	264	&	45	&&	4.11E-07	&	4.12E-07	&	4.55E-07	\\
	&	6	&	1041	&	76	&	50	&&	4.11E-07	&	4.11E-07	&	4.54E-07	\\
	&	7	&	1041	&	256	&	55	&&	4.10E-07	&	4.11E-07	&	4.54E-07	\\
	&	8	&	1041	&	270	&	80	&&	4.04E-07	&	4.11E-07	&	4.53E-07	\\
	&	9	&	1041	&	129	&	84	&&	4.04E-07	&	4.11E-07	&	4.53E-07	\\
	&	10	&	1041	&	218	&	87	&&	4.04E-07	&	4.11E-07	&	4.53E-07	\\
\midrule
2	&	1	&	1041	&	31	&	14	&&	2.43E-07	&	2.43E-07	&	2.43E-07	\\
	&	2	&	1181	&	93	&	13	&&	2.39E-07	&	2.38E-07	&	2.43E-07	\\
	&	3	&	1041	&	106	&	36	&&	2.36E-07	&	2.37E-07	&	2.42E-07	\\
	&	4	&	1041	&	289	&	33	&&	2.36E-07	&	2.37E-07	&	2.42E-07	\\
	&	5	&	1041	&	232	&	46	&&	2.36E-07	&	2.37E-07	&	2.42E-07	\\
	&	6	&	1041	&	141	&	56	&&	2.36E-07	&	2.37E-07	&	2.41E-07	\\
	&	7	&	1041	&	106	&	71	&&	2.36E-07	&	2.37E-07	&	2.41E-07	\\
	&	8	&	1041	&	226	&	69	&&	2.36E-07	&	2.37E-07	&	2.41E-07	\\
	&	9	&	1041	&	565	&	97	&&	2.36E-07	&	2.37E-07	&	2.41E-07	\\
	&	10	&	1041	&	260	&	92	&&	2.36E-07	&	2.37E-07	&	2.41E-07	\\
\midrule
3	&	1	&	1041	&	37	&	15	&&	4.99E-07	&	4.99E-07	&	4.99E-07	\\
	&	2	&	1381	&	68	&	19	&&	3.21E-07	&	3.27E-07	&	4.97E-07	\\
	&	3	&	1041	&	64	&	35	&&	3.21E-07	&	3.18E-07	&	4.95E-07	\\
	&	4	&	1041	&	121	&	32	&&	3.21E-07	&	3.15E-07	&	4.94E-07	\\
	&	5	&	1041	&	223	&	73	&&	3.21E-07	&	3.15E-07	&	4.93E-07	\\
	&	6	&	1041	&	92	&	66	&&	3.21E-07	&	3.15E-07	&	4.91E-07	\\
	&	7	&	1041	&	192	&	92	&&	3.21E-07	&	3.15E-07	&	4.90E-07	\\
	&	8	&	1041	&	1430	&	53	&&	3.21E-07	&	3.13E-07	&	4.90E-07	\\
	&	9	&	1041	&	1232	&	942	&&	3.21E-07	&	3.13E-07	&	3.43E-07	\\
	&	10	&	1041	&	451	&	324	&&	3.21E-07	&	3.13E-07	&	3.40E-07	\\
\midrule
4	&	1	&	1041	&	32	&	26	&&	3.86E-06	&	3.86E-06	&	3.86E-06	\\
	&	2	&	1041	&	73	&	100	&&	2.74E-06	&	2.74E-06	&	2.74E-06	\\
	&	3	&	1041	&	89	&	23	&&	2.74E-06	&	2.66E-06	&	2.74E-06	\\
	&	4	&	1241	&	208	&	38	&&	2.68E-06	&	2.63E-06	&	2.74E-06	\\
	&	5	&	1041	&	70	&	30	&&	2.67E-06	&	2.62E-06	&	2.74E-06	\\
	&	6	&	1041	&	157	&	34	&&	2.67E-06	&	2.61E-06	&	2.74E-06	\\
	&	7	&	1041	&	92	&	30	&&	2.67E-06	&	2.61E-06	&	2.74E-06	\\
	&	8	&	1041	&	219	&	44	&&	2.67E-06	&	2.61E-06	&	2.74E-06	\\
	&	9	&	1041	&	108	&	63	&&	2.61E-06	&	2.61E-06	&	2.74E-06	\\
	&	10	&	1041	&	308	&	50	&&	2.61E-06	&	2.61E-06	&	2.74E-06	\\
\midrule
5	&	1	&	1041	&	19	&	18	&&	1.93E-06	&	1.93E-06	&	1.93E-06	\\
	&	2	&	1041	&	80	&	96	&&	1.24E-06	&	1.24E-06	&	1.24E-06	\\
	&	3	&	1141	&	171	&	33	&&	1.22E-06	&	1.22E-06	&	1.24E-06	\\
	&	4	&	1161	&	210	&	30	&&	1.22E-06	&	1.21E-06	&	1.24E-06	\\
	&	5	&	1041	&	136	&	43	&&	1.22E-06	&	1.21E-06	&	1.24E-06	\\
	&	6	&	1041	&	280	&	34	&&	1.21E-06	&	1.21E-06	&	1.24E-06	\\
	&	7	&	1041	&	285	&	68	&&	1.21E-06	&	1.21E-06	&	1.24E-06	\\
	&	8	&	1041	&	338	&	58	&&	1.21E-06	&	1.21E-06	&	1.24E-06	\\
	&	9	&	1041	&	103	&	80	&&	1.21E-06	&	1.21E-06	&	1.24E-06	\\
	&	10	&	1041	&	169	&	92	&&	1.21E-06	&	1.21E-06	&	1.24E-06	\\
    \bottomrule
    \end{tabular}%
  \label{table:dataMatIIwarm}%
\end{table}%

\begin{table}[htbp]\scriptsize
  \centering
  \caption{Number of function calls and objective values found for Material Set III for algorithms GA$_r$, MADS$_r$, and RAGS$_r$.}
    \begin{tabular}{@{}cccccccccc@{}}
    \toprule
     && \multicolumn{3}{c}{Number of function calls}  && \multicolumn{3}{c}{Objective value}  \\[-0.2mm]
    \midrule
    Case  & nd    & GA    & MADS  & RAGS  &       & GA    & MADS  & RAGS \\[-0.2mm]
    \midrule
    1	&	1	&	385	&	293	&	159	&	&	5.39E-07	&	5.42E-07	&	5.37E-07	\\
	&	2	&	1015	&	1158	&	635	&	&	5.22E-07	&	5.16E-07	&	5.16E-07	\\
	&	3	&	1855	&	2408	&	1631	&	&	5.03E-07	&	4.81E-07	&	4.79E-07	\\
	&	4	&	2835	&	3891	&	2823	&	&	5.04E-07	&	4.80E-07	&	4.68E-07	\\
	&	5	&	3885	&	5849	&	4165	&	&	4.84E-07	&	4.74E-07	&	4.60E-07	\\
	&	6	&	4935	&	8353	&	5533	&	&	4.83E-07	&	4.61E-07	&	4.62E-07	\\
	&	7	&	5915	&	10105	&	7017	&	&	4.84E-07	&	4.68E-07	&	4.59E-07	\\
	&	8	&	6755	&	12987	&	7926	&	&	4.70E-07	&	4.62E-07	&	4.63E-07	\\
	&	9	&	7385	&	14166	&	9089	&	&	4.72E-07	&	4.68E-07	&	4.59E-07	\\
	&	10	&	350	&	1365	&	610	&	&	4.87E-07	&	4.65E-07	&	4.59E-07	\\
\midrule																
2	&	1	&	395	&	273	&	117	&	&	4.70E-07	&	4.64E-07	&	4.62E-07	\\
	&	2	&	1055	&	1127	&	302	&	&	4.41E-07	&	4.40E-07	&	4.54E-07	\\
	&	3	&	1910	&	2559	&	836	&	&	4.39E-07	&	4.40E-07	&	4.50E-07	\\
	&	4	&	2890	&	4558	&	1562	&	&	4.40E-07	&	4.37E-07	&	4.43E-07	\\
	&	5	&	3940	&	6421	&	3539	&	&	4.38E-07	&	4.35E-07	&	4.40E-07	\\
	&	6	&	4990	&	8294	&	5246	&	&	4.38E-07	&	4.36E-07	&	4.39E-07	\\
	&	7	&	5970	&	10386	&	6726	&	&	4.38E-07	&	4.38E-07	&	4.39E-07	\\
	&	8	&	6810	&	13105	&	8304	&	&	4.40E-07	&	4.40E-07	&	4.38E-07	\\
	&	9	&	7440	&	15972	&	9417	&	&	4.40E-07	&	4.41E-07	&	4.38E-07	\\
	&	10	&	350	&	1038	&	715	&	&	4.50E-07	&	4.51E-07	&	4.38E-07	\\
\midrule																
3	&	1	&	415	&	290	&	156	&	&	8.98E-07	&	8.41E-07	&	8.41E-07	\\
	&	2	&	1075	&	946	&	702	&	&	8.09E-07	&	4.97E-07	&	4.91E-07	\\
	&	3	&	1945	&	1930	&	1494	&	&	6.49E-07	&	4.94E-07	&	4.91E-07	\\
	&	4	&	2985	&	3352	&	2573	&	&	5.04E-07	&	4.93E-07	&	4.91E-07	\\
	&	5	&	4110	&	5335	&	3581	&	&	5.28E-07	&	4.85E-07	&	5.00E-07	\\
	&	6	&	5220	&	7326	&	5248	&	&	4.97E-07	&	4.95E-07	&	5.08E-07	\\
	&	7	&	6270	&	8942	&	7197	&	&	5.11E-07	&	4.85E-07	&	5.09E-07	\\
	&	8	&	7190	&	11827	&	9596	&	&	5.04E-07	&	4.90E-07	&	5.10E-07	\\
	&	9	&	7820	&	13554	&	13366	&	&	5.14E-07	&	4.92E-07	&	5.09E-07	\\
	&	10	&	350	&	1306	&	955	&	&	5.78E-07	&	4.88E-07	&	5.11E-07	\\
\midrule																
4	&	1	&	425	&	317	&	202	&	&	8.63E-06	&	8.64E-06	&	8.62E-06	\\
	&	2	&	1075	&	989	&	1064	&	&	8.65E-06	&	6.51E-06	&	8.63E-06	\\
	&	3	&	1945	&	2143	&	2196	&	&	6.11E-06	&	6.10E-06	&	8.67E-06	\\
	&	4	&	2965	&	4028	&	3572	&	&	6.17E-06	&	5.87E-06	&	8.70E-06	\\
	&	5	&	4015	&	6565	&	5950	&	&	6.32E-06	&	5.84E-06	&	8.70E-06	\\
	&	6	&	5065	&	10071	&	8418	&	&	8.73E-06	&	5.80E-06	&	8.73E-06	\\
	&	7	&	6150	&	12257	&	10836	&	&	6.47E-06	&	5.81E-06	&	8.77E-06	\\
	&	8	&	7030	&	15980	&	12886	&	&	6.35E-06	&	5.79E-06	&	8.78E-06	\\
	&	9	&	7930	&	18156	&	14585	&	&	6.36E-06	&	5.72E-06	&	8.81E-06	\\
	&	10	&	500	&	1194	&	1098	&	&	6.19E-06	&	5.71E-06	&	8.82E-06	\\
\midrule																
5	&	1	&	465	&	298	&	284	&	&	6.85E-06	&	6.21E-06	&	6.19E-06	\\
	&	2	&	1125	&	1163	&	813	&	&	5.60E-06	&	3.20E-06	&	3.20E-06	\\
	&	3	&	2070	&	2332	&	1586	&	&	3.60E-06	&	2.46E-06	&	2.47E-06	\\
	&	4	&	3110	&	3957	&	2438	&	&	3.40E-06	&	2.41E-06	&	2.43E-06	\\
	&	5	&	4310	&	6951	&	3373	&	&	2.68E-06	&	2.37E-06	&	2.41E-06	\\
	&	6	&	5420	&	9313	&	4481	&	&	2.56E-06	&	2.36E-06	&	2.39E-06	\\
	&	7	&	6470	&	12817	&	5570	&	&	2.41E-06	&	2.36E-06	&	2.38E-06	\\
	&	8	&	7350	&	16541	&	6600	&	&	2.41E-06	&	2.35E-06	&	2.38E-06	\\
	&	9	&	7980	&	19290	&	7276	&	&	2.57E-06	&	2.35E-06	&	2.38E-06	\\
	&	10	&	350	&	2299	&	367	&	&	2.66E-06	&	2.34E-06	&	2.39E-06	\\
    \bottomrule
    \end{tabular}%
  \label{table:dataMatIIIrandom}%
\end{table}%

\begin{table}[h!]\scriptsize
  \centering
  \caption{Number of function calls and objective values found for Material Set III for algorithms GA$_w$, MADS$_w$, and RAGS$_w$.}
    \begin{tabular}{@{}cccccccccc@{}}
    \toprule
     && \multicolumn{3}{c}{Number of function calls}  && \multicolumn{3}{c}{Objective value}  \\[-0.2mm]
    \midrule    Case  & nd    & GA$_w$    & MADS$_w$  & RAGS$_w$  &       & GA$_w$    & MADS$_w$  & RAGS$_w$ \\[-0.2mm]
    \midrule
1	&	1	&	1041	&	36	&	29	&&	5.37E-07	&	5.42E-07	&	5.37E-07	\\
	&	2	&	1041	&	84	&	37	&&	5.16E-07	&	5.17E-07	&	5.34E-07	\\
	&	3	&	1961	&	141	&	47	&&	4.89E-07	&	4.80E-07	&	5.32E-07	\\
	&	4	&	1041	&	205	&	501	&&	4.88E-07	&	4.74E-07	&	5.17E-07	\\
	&	5	&	1041	&	187	&	72	&&	4.88E-07	&	4.69E-07	&	5.16E-07	\\
	&	6	&	1041	&	103	&	43	&&	4.87E-07	&	4.69E-07	&	5.16E-07	\\
	&	7	&	1041	&	441	&	38	&&	4.87E-07	&	4.68E-07	&	5.15E-07	\\
	&	8	&	1041	&	162	&	85	&&	4.87E-07	&	4.68E-07	&	5.15E-07	\\
	&	9	&	1041	&	584	&	75	&&	4.86E-07	&	4.68E-07	&	5.15E-07	\\
	&	10	&	1041	&	839	&	59	&&	4.86E-07	&	4.68E-07	&	5.15E-07	\\
\midrule
2	&	1	&	1041	&	29	&	15	&&	4.62E-07	&	4.64E-07	&	4.62E-07	\\
	&	2	&	1081	&	69	&	32	&&	4.39E-07	&	4.40E-07	&	4.62E-07	\\
	&	3	&	1041	&	83	&	27	&&	4.39E-07	&	4.35E-07	&	4.62E-07	\\
	&	4	&	1041	&	68	&	38	&&	4.39E-07	&	4.35E-07	&	4.61E-07	\\
	&	5	&	1041	&	108	&	48	&&	4.39E-07	&	4.35E-07	&	4.61E-07	\\
	&	6	&	1041	&	131	&	80	&&	4.38E-07	&	4.34E-07	&	4.61E-07	\\
	&	7	&	1041	&	185	&	68	&&	4.39E-07	&	4.34E-07	&	4.61E-07	\\
	&	8	&	1041	&	178	&	78	&&	4.39E-07	&	4.34E-07	&	4.60E-07	\\
	&	9	&	1041	&	523	&	71	&&	4.39E-07	&	4.34E-07	&	4.60E-07	\\
	&	10	&	1041	&	152	&	96	&&	4.39E-07	&	4.34E-07	&	4.60E-07	\\
\midrule
3	&	1	&	1041	&	29	&	12	&&	8.41E-07	&	8.41E-07	&	8.41E-07	\\
	&	2	&	1041	&	53	&	30	&&	5.01E-07	&	4.94E-07	&	8.39E-07	\\
	&	3	&	1041	&	65	&	24	&&	5.01E-07	&	4.93E-07	&	8.37E-07	\\
	&	4	&	1041	&	56	&	41	&&	4.95E-07	&	4.93E-07	&	8.36E-07	\\
	&	5	&	1041	&	90	&	63	&&	4.88E-07	&	4.93E-07	&	8.34E-07	\\
	&	6	&	1041	&	55	&	48	&&	4.87E-07	&	4.93E-07	&	8.33E-07	\\
	&	7	&	1041	&	300	&	57	&&	4.87E-07	&	4.92E-07	&	8.33E-07	\\
	&	8	&	1041	&	199	&	59	&&	4.87E-07	&	4.91E-07	&	8.32E-07	\\
	&	9	&	1041	&	416	&	123	&&	4.87E-07	&	4.91E-07	&	8.32E-07	\\
	&	10	&	1041	&	291	&	68	&&	4.87E-07	&	4.91E-07	&	8.32E-07	\\
\midrule
4	&	1	&	1041	&	37	&	19	&&	8.62E-06	&	8.64E-06	&	8.62E-06	\\
	&	2	&	1041	&	54	&	18	&&	8.62E-06	&	8.62E-06	&	8.62E-06	\\
	&	3	&	1041	&	63	&	33	&&	8.62E-06	&	8.62E-06	&	8.62E-06	\\
	&	4	&	1041	&	59	&	42	&&	8.62E-06	&	8.62E-06	&	8.63E-06	\\
	&	5	&	1041	&	118	&	64	&&	8.62E-06	&	8.62E-06	&	8.64E-06	\\
	&	6	&	1041	&	319	&	51	&&	8.62E-06	&	8.62E-06	&	8.66E-06	\\
	&	7	&	1041	&	41	&	71	&&	8.62E-06	&	8.62E-06	&	8.68E-06	\\
	&	8	&	1041	&	1360	&	62	&&	8.62E-06	&	5.92E-06	&	8.71E-06	\\
	&	9	&	1041	&	1008	&	75	&&	8.62E-06	&	5.72E-06	&	8.73E-06	\\
	&	10	&	1041	&	753	&	81	&&	8.62E-06	&	5.68E-06	&	8.76E-06	\\
\midrule
5	&	1	&	1041	&	26	&	30	&&	6.19E-06	&	6.21E-06	&	6.19E-06	\\
	&	2	&	1041	&	74	&	24	&&	3.20E-06	&	3.20E-06	&	6.17E-06	\\
	&	3	&	1041	&	160	&	31	&&	2.46E-06	&	2.46E-06	&	6.15E-06	\\
	&	4	&	1201	&	162	&	20	&&	2.42E-06	&	2.42E-06	&	6.14E-06	\\
	&	5	&	1261	&	286	&	630	&&	2.37E-06	&	2.37E-06	&	2.49E-06	\\
	&	6	&	1041	&	281	&	321	&&	2.37E-06	&	2.36E-06	&	2.42E-06	\\
	&	7	&	1061	&	456	&	68	&&	2.35E-06	&	2.36E-06	&	2.42E-06	\\
	&	8	&	1041	&	459	&	59	&&	2.34E-06	&	2.36E-06	&	2.42E-06	\\
	&	9	&	1041	&	204	&	56	&&	2.34E-06	&	2.35E-06	&	2.42E-06	\\
	&	10	&	1041	&	640	&	106	&&	2.34E-06	&	2.34E-06	&	2.42E-06	\\
    \bottomrule
    \end{tabular}%
  \label{table:dataMatIIIwarm}%
\end{table}%

\end{document}